\documentclass[a4paper]{article}
\usepackage{mathpazo}
\usepackage[english]{babel}
\usepackage[utf8]{inputenc}
\usepackage[T1]{fontenc}
\usepackage{amsfonts,amsmath,amssymb,amsthm}
\usepackage{graphicx}
\usepackage{color}
\usepackage{hyperref}
\usepackage{enumitem}
\usepackage{tikz}
\usepackage{blkarray}
\usepackage{array}
\usepackage{algorithm}
\usepackage{algpseudocode}
\usetikzlibrary{decorations.markings}

\setlist[enumerate]{itemsep=0.1cm}

\newtheorem{theorem}{Theorem}[section]

\newtheorem{proposition}[theorem]{Proposition}
\newtheorem{lemma}[theorem]{Lemma}

\newtheorem*{corollary*}{Corollary}

\theoremstyle{definition}
\newtheorem{definition}{Definition}[section]
\newtheorem{example}{Example}[section]

\theoremstyle{remark}
\newtheorem*{remark}{Remark}

\newcommand{\Z}{\mathbf Z}

\newcommand{\GG}{\mathcal G}

\newcommand{\ra}{\rightarrow}

\newcommand{\dfl}[1]{\text{dFl}(#1)}

\title{Simplicial \(q\)-connectivity of directed graphs with applications to network analysis}
\author{Henri Riihim\"{a}ki\thanks{Department of Mathematics, KTH Royal Institute of Technology, \url{henrir@kth.se}}}
\begin{document}

\maketitle

\begin{abstract}
	Directed graphs are ubiquitous models for networks, and topological spaces they generate, such as the directed flag complex, have become useful objects in applied topology. The simplices are formed from directed cliques. We extend Atkin's theory of \(q\)-connectivity to the case of directed simplices. This results in a preorder where simplices are related by sequences of simplices that share a \(q\)-face with respect to directions specified by chosen face maps.  We leverage the Alexandroff equivalence between preorders and  topological spaces to introduce a new class of topological spaces for directed graphs, enabling to assign new homotopy types different from those of directed flag complexes as seen by simplicial homology. We further introduce simplicial path analysis enabled by the connectivity preorders. As an application we characterise structural differences between various brain networks by computing their longest simplicial paths.
\end{abstract}

\section{Introduction}
Directed graphs, which we usually call \emph{digraphs} in this paper, are a model for various phenomena in the sciences, for example networks of neurons in the brain or gene regulatory networks. In topological data analysis (TDA) particularly the advent of applying topological tools to questions in neuroscience has spawned interest in constructing topological spaces out of digraphs, developing computational tools for obtaining topological information, and using these to understand networks and phenomena they support. For a progression of works on these ideas, see \cite{Tribes_math,Flagser_paper,Frontiers_paper,Tribes_neuroscience}. Our main example of a topological space on a digraph \(\GG\) is the \emph{directed flag complex}, which is constructed from the directed cliques of \(\GG\). For example, a 2-simplex is given by an ordered sequence of vertices \((v_0,v_1,v_2)\) whenever any ordered pair \((v_i,v_j)\), for \(i < j\), is a directed edge in \(\GG\). By construction the simplices are endowed with an inherent directionality. For a recent work on computing the homotopy type of the directed flag complex of the C. elegans neuronal network, see \cite{Govc_2020}.

\begin{figure}[h]
	\begin{center}
		\begin{tikzpicture}
		\tikzstyle{point}=[circle,thick,draw=black,fill=black,inner sep=0pt,minimum width=2pt,minimum height=2pt]
		\tikzset{dirarrow/.style={postaction={decorate,decoration={markings,mark=at position .7 with {\arrow[#1]{stealth}}}}}}
		\tikzstyle{arc}=[shorten >= 8pt,shorten <= 8pt,->]
		
		\coordinate (a) at (-0.25,0);
		\coordinate (b) at (1,0.5);
		\coordinate (c) at (0.75,-0.5);
		\coordinate (d) at (-0.75,-1);
		\coordinate (e) at (2.25,-0.25);
		\coordinate (f) at (2.5,0.75);
		
		\draw[fill=magenta,opacity=0.3] (a) -- (b) -- (c) --cycle; 
		
		\draw[fill=magenta,opacity=0.3] (a) -- (d) -- (c) --cycle; 
		
		\draw[fill=magenta,opacity=0.3] (b) -- (c) -- (e) --cycle;
		
		\draw[fill=magenta,opacity=0.3] (b) -- (f) -- (e) --cycle;
		
		\draw[postaction={dirarrow}] (a) to (b);
		\draw[postaction={dirarrow}] (a) -- (c);
		\draw[postaction={dirarrow}] (d) -- (a);
		\draw[postaction={dirarrow}] (d) -- (c);
		\draw[postaction={dirarrow}] (c) -- (b);
		\draw[postaction={dirarrow}] (b) -- (e);
		\draw[postaction={dirarrow}] (c) -- (e);
		\draw[postaction={dirarrow}] (b) -- (f);
		\draw[postaction={dirarrow}] (e) -- (f);
		
		\node[point] at (a) {};
		\node[point] at (d) {};
		\node[point] at (c) {};
		\node[point] at (b) {};
		\node[point] at (e) {};
		\node[point] at (f) {};
		
		\draw[arc,blue] (d) .. controls (0.5,0.3) and (1.5,-0.2) .. (f);
		\coordinate (a') at (4.75,0);
		\coordinate (b') at (6,0.5);
		\coordinate (c') at (5.75,-0.5);
		\coordinate (d') at (4.25,-1);
		\coordinate (e') at (7.25,-0.25);
		\coordinate (f') at (7.5,0.75);
		
		\draw[fill=magenta,opacity=0.3] (a') -- (b') -- (c') --cycle; 
		
		\draw[fill=magenta,opacity=0.3] (a') -- (d') -- (c') --cycle; 
		
		\draw[fill=magenta,opacity=0.3] (b') -- (c') -- (e') --cycle;
		
		\draw[fill=magenta,opacity=0.3] (b') -- (f') -- (e') --cycle;
		
		\draw[postaction={dirarrow}] (a') to (b');
		\draw[postaction={dirarrow}] (a') -- (c');
		\draw[postaction={dirarrow}] (a') -- (d');
		\draw[postaction={dirarrow}] (c') -- (d');
		\draw[postaction={dirarrow}] (c') -- (b');
		\draw[postaction={dirarrow}] (b') -- (e');
		\draw[postaction={dirarrow}] (c') -- (e');
		\draw[postaction={dirarrow}] (f') -- (b');
		\draw[postaction={dirarrow}] (f') -- (e');
		
		\node[point] at (a') {};
		\node[point] at (d') {};
		\node[point] at (c') {};
		\node[point] at (b') {};
		\node[point] at (e') {};
		\node[point] at (f') {};
		
		\draw[arc,blue] (a') .. controls (6,0.1) and (6,0) .. (e');
		\draw[arc,blue] (a') .. controls (5,-0.5) .. (d');
		\draw[arc,blue] (f') .. controls (6.9,0.4) .. (e');
		\end{tikzpicture}
	\end{center}
	\caption{The directed simplices on the left complex are equidirected such that they generate a path of 2-simplices. On the right the directionality of simplices acts as an obstruction and thus creates multiple shorter paths.}
	\label{fig:example_pathways}
\end{figure}
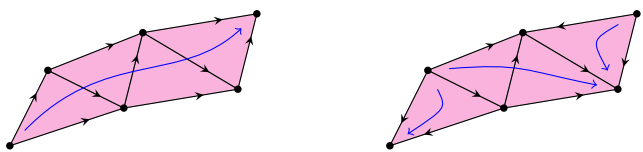

This paper stems from two streams of ideas trying to answer the questions below, relating to further topological and combinatorial analyses of directed flag complexes of digraphs and networks:
\begin{enumerate}
	\item[Q1] We can construct the directed cliques from digraphs as sets. What other topological spaces can we construct from these sets, that would be different from the directed flag complex and its homotopy type, and would be more refined constructions in terms of the underlying digraphs? We will show in Example \ref{ex:2_spheres} two different digraphs whose directed flag complexes are topologically 2-spheres. This stems from the fact that even though the directed flag complex is constructed from directed cliques, the topological space is that of the associated geometric realisation, and hence the directionality information is lost.
	\item[Q2] Digraphs naturally support various directed edge paths. But since directed cliques have an inherent directionality, can we extend to \emph{higher simplicial paths} formed by sequences of directed simplices? Figure \ref{fig:example_pathways} gives a visual illustration on this. The 2-dimensional complexes in the figure are also further examples of two topologically indistinguishable, in fact contractible, spaces; the simplicial directionality is nonetheless different as the arrows illustrate.
\end{enumerate}

A convenient framework to answer these questions comes from so-called \emph{\(q\)-connectivity}. In a series of works in the 1970's Atkin defined and studied what he called \(Q\)-analysis of simplicial complexes \cite{Atkin1972,Atkin1974,Atkin1974_book,Atkin1976}. In \(Q\)-analysis simplices are \(q\)-near if they share a common \(q\)-dimensional face. This relation is transitively closed to define the \(q\)-connectivity equivalence relation by sequences of \(q\)-near simplices. The equivalence classes capture \(q\)-connected components of the simplicial complex, a direct generalisation of path components. The approach of \(Q\)-analysis has been used in various applications such as social sciences \cite{Atkin1972,Atkin1974_book}, network analysis \cite{Laubenbacher_networks} and analysis of neuroscientific data \cite{Tadic_et_al}, and it has given rise to a combinatorial homotopy theory \cite{Barcelo_etal_2001,BarceloLaubenbacher_2005,KramerLaubenbacher1998}.

Adapting \(Q\)-analysis to the setting of directed simplices, as is the case with directed flag complexes, is the main content of our paper:
\begin{enumerate}
	\item[Q3] Classical \(Q\)-analysis works in the setting of simplicial complexes. Hence it can fail to distinguish very different directed flag complexes, see Example \ref{ex:q_connectivity_problem}. What is the appropriate extension of \(q\)-connectivity to directed simplices, allowing one to apply directed notion of \(Q\)-analysis to digraphs?
\end{enumerate}

Our directed notion of \(q\)-nearness is based on two directed simplices \(\sigma\) and \(\tau\) sharing a \(q\)-dimensional face, but we also impose a condition on \emph{how} the face is shared with respect to the directionality of \(\sigma\) and \(\tau\). We use the simplicial face maps \(d_i\) to indicate a face in direction \(i\): if \(\sigma\) is an \(n\)-simplex \((v_0,\dots,v_n)\), then \(d_i(\sigma) = (v_0,\dots,\widehat{v_i},\dots, v_n)\) is an \((n-1)\)-face where \(\widehat{v_i}\) denotes the removal of the vertex \(v_i\). To accommodate \(q\)-connectivity we define a slightly modified version of the face maps, \(\widehat{d_i}\), in Definition \ref{def:new_face_map}. Our directed \(q\)-nearness between an ordered pair of simplices \((\sigma,\tau)\) is now denoted by a triple \((q,\widehat{d_i},\widehat{d_j})\), designating that \(\widehat{d_i}(\sigma) \hookleftarrow \alpha \hookrightarrow \widehat{d_j}(\tau)\) for some \(q\)-simplex \(\alpha\), i.e.\ \(\sigma\) and \(\tau\) share a \(q\)-face in the directions specified by the new face maps \(\widehat{d_i}\) and \(\widehat{d_j}\) (Definition \ref{def:q_near_directed}). This answers Question 3 and Section \ref{sec:directed_q_analysis} covers the theory of directed \(Q\)-analysis.

Closing the directed \(q\)-nearness transitively gives us directed \(q\)-connectivity. The relation obtained is a preorder between simplices. This is the crucial difference from classical \(Q\)-analysis. Instead of identifying \(q\)-connected components of the equivalence relation, we can regard the preorder of directed \(q\)-connectivity as a digraph of simplices. Question 2 above is then answered by path searches in this digraph. To answer Question 1 we leverage the theory of finite topological spaces \cite{Barmak_book}. By a classical result of Alexandroff, preorders are in bijection with topological spaces with Alexandroff topologies. Applying this to the \(q\)-connectivity preorders allows us to assign new homotopy types to directed flag complexes.

We hope that the analysis of simplicial paths might become a new and interesting tool for network science. As such we view this as a form of combinatorial dimensionality reduction: instead of the path structure of the digraph we look at the path structure of higher dimensional simplices in the \(q\)-connectivity preorder. For example, on the left complex of Figure \ref{fig:example_pathways} there are 8 edge paths from the leftmost vertex to the rightmost vertex, compared to the single path of directed 2-simplices indicated by the blue arrow. Finding motifs and other hierarchical structures is a prominent theme in network science and in its applications to neuroscience, see for example the survey \cite{Curto_motif}. The simplicial paths might thus provide a new higher-dimensional motif. As our main application to network analysis in Section \ref{sec:applications} we take a step in this direction by characterising various neuronal networks in terms of their longest simplicial paths with respect to different directed \(q\)-connectivities. Recently a cascading dynamics on networks was studied on the level of simplicial complexes, and it was shown that this dynamics follows geometric channels of simplices \cite{simplicial_cascades}. The simplicial paths are an appropriate framework for such an analysis in the case of directed networks.

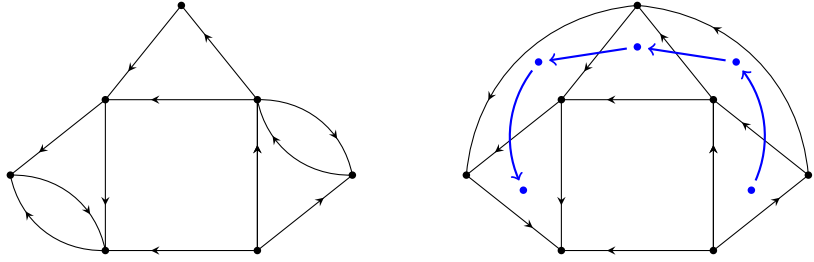
\begin{figure}[h]
	\begin{center}
		\begin{tikzpicture}
		\tikzstyle{point}=[circle,thick,draw=black,fill=black,inner sep=0pt,minimum width=2pt,minimum height=2pt]
		\tikzstyle{bluepoint}=[circle,thick,draw=blue,fill=blue,inner sep=0pt,minimum width=2pt,minimum height=2pt]
		\tikzset{dirarrow/.style={postaction={decorate,decoration={markings,mark=at position .7 with {\arrow[#1]{stealth}}}}}}
		\tikzstyle{arc}=[shorten >= 4pt,shorten <= 4pt,->]
		
		\coordinate (a) at (0,0);
		\coordinate (b) at (-2,0);
		\coordinate (c) at (-3.25,1);
		\coordinate (d) at (-2,2);
		\coordinate (e) at (0,2);
		\coordinate (f) at (1.25,1);
		\coordinate (g) at (-1,3.25);
		
		
		\draw[postaction={dirarrow}] (a) -- (b);
		\draw[postaction={dirarrow}] (a) -- (e);
		\draw[postaction={dirarrow}] (b) to [bend left=40] (c);
		\draw[postaction={dirarrow}] (c) to [bend left=40] (b);
		\draw[postaction={dirarrow}] (d) -- (c);
		\draw[postaction={dirarrow}] (d) -- (b);
		\draw[postaction={dirarrow}] (a) -- (e);
		\draw[postaction={dirarrow}] (a) -- (f);
		\draw[postaction={dirarrow}] (f) to [bend left=40] (e);
		\draw[postaction={dirarrow}] (e) to [bend left=40] (f);
		\draw[postaction={dirarrow}] (e) -- (d);
		\draw[postaction={dirarrow}] (g) -- (d);
		\draw[postaction={dirarrow}] (e) -- (g);
		
		\node[point] at (a) {};
		\node[point] at (d) {};
		\node[point] at (c) {};
		\node[point] at (b) {};
		\node[point] at (e) {};
		\node[point] at (f) {};
		\node[point] at (g) {};
		\coordinate (a') at (6,0);
		\coordinate (b') at (4,0);
		\coordinate (c') at (2.75,1);
		\coordinate (d') at (4,2);
		\coordinate (e') at (6,2);
		\coordinate (f') at (7.25,1);
		\coordinate (g') at (5,3.25);
		
		
		\draw[postaction={dirarrow}] (a') -- (b');
		\draw[postaction={dirarrow}] (a') -- (e');
		\draw[postaction={dirarrow}] (c') -- (b');
		\draw[postaction={dirarrow}] (d') -- (c');
		\draw[postaction={dirarrow}] (d') -- (b');
		\draw[postaction={dirarrow}] (a') -- (e');
		\draw[postaction={dirarrow}] (a') -- (f');
		\draw[postaction={dirarrow}] (f') -- (e');
		\draw[postaction={dirarrow}] (f') to [bend right=40] (g');
		\draw[postaction={dirarrow}] (e') -- (d');
		\draw[postaction={dirarrow}] (g') -- (d');
		\draw[postaction={dirarrow}] (e') -- (g');
		\draw[postaction={dirarrow}] (g') to [bend right=40] (c');
		
		\node[point] at (a') {};
		\node[point] at (d') {};
		\node[point] at (c') {};
		\node[point] at (b') {};
		\node[point] at (e') {};
		\node[point] at (f') {};
		\node[point] at (g') {};
		
		\coordinate (p1) at (6.5,0.8);
		\coordinate (p2) at (6.3,2.5);
		\coordinate (p3) at (5,2.7);
		\coordinate (p4) at (3.7,2.5);
		\coordinate (p5) at (3.5,0.8);
		
		\node[bluepoint] at (p1) {};
		\node[bluepoint] at (p2) {};
		\node[bluepoint] at (p3) {};
		\node[bluepoint] at (p4) {};
		\node[bluepoint] at (p5) {};
		
		\draw[blue,thick,arc] (p1) to [bend right=30] (p2);
		\draw[blue,thick,arc] (p2) -- (p3);
		\draw[blue,thick,arc] (p3) -- (p4);
		\draw[blue,thick,arc] (p4) to [bend right=30] (p5);
		\end{tikzpicture}
	\end{center}
	\caption{Two digraphs with the same simplex counts and homology, but on the right we see a simplicial path as drawn in blue arrows.}
	\label{fig:example_path_analysis_same_homology}
\end{figure}

To further explain why simplicial path analysis of digraphs might be an interesting tool, complimenting more standard tools such as homology computations of directed flag complexes, we again resort to visual aids in Figure \ref{fig:example_path_analysis_same_homology}. The directed flag complexes of the digraphs shown have the same simplex counts (7 vertices, 12 edges, 5 2-simplices) and the same homology with \(\Z_2\) coefficients (Betti numbers 1 and 1 in degrees 0 and 1, respectively). The combinatorics of the 2-simplices, however, is very different. The right digraph has in our parlance a \((1,\widehat{d_0},\widehat{d_2})\)-connected path of 2-simplices as depicted with blue arrows.

The approach offered by \(q\)-connectivity can be conceptually understood as follows. The face poset of a simplicial complex has as elements the simplices and the partial order relation is given by simplicial face inclusions. It is a particular feature of face posets that they organise into levels determined by the dimensions of the simplices, with order relations only going from lower levels to higher levels. What \(q\)-connectivity then does is to add "horizontal" relations between simplices provided that they share a face in a lower level. This is made explicit in Sections \ref{sec:q_analysis} and \ref{sec:directed_q_analysis} by showing that the face poset is a subrelation of \(q\)-connectivity, both in undirected and directed cases. It is well known that the homotopy types of a simplicial complex and its face poset agree \cite{Barmak_book}. The modification of the latter by \(q\)-connectivity then conceptually explains why it can extract new kind of topological information able to distinguish spaces that, for example, simplicial homology cannot as remarked in Question 1 above.

This paper is organised as follows. Section \ref{sec:graphs_complexes} reviews basic definitions of digraphs and simplicial complexes, and introduces the directed flag complex. Section \ref{sec:q_analysis} gives a self-contained introduction to Atkin's classical \(Q\)-analysis. We give special emphasis to pseudomanifolds and prove some characterisations in terms of connectivity graphs arising from \(Q\)-analysis. In this we are motivated by the neural manifold hypothesis \cite{neural_manifolds} which assumes that neural activity data of a collection of neurons resides on a lower-dimensional parametrising manifold. In this context it might be of interest to ask whether a collection of active neurons constitutes a combinatorial pseudomanifold on the network level. In Section \ref{sec:q_analysis} we also show the relationship between face poset and \(q\)-connectivity, as well as the observation that \(k\)-clique communities known from network analysis are subsumed by the latter. Section \ref{sec:directed_q_analysis} develops the directed \(q\)-connectivity and proves its basic properties. In Section \ref{sec:dir_struct_sys} we outline the approach from finite (Alexandroff) spaces to study \(q\)-connectivities topologically. Again we give emphasis to the associated connectivity digraphs. We finish with Section \ref{sec:applications} where we apply the new directed \(q\)-connectivity and simplicial path analysis to various brain networks. It is shown that different networks exhibit quite different connectivity structures which can serve as new structural fingerprints for subsequent analysis, such as featurisations for machine learning tasks.

\section{Graphs and complexes}\label{sec:graphs_complexes}
We start by fixing the notions related to graphs and simplicial complexes. In this paper we are interested in finite graphs and finite complexes.

\begin{definition}\label{def:graph}
	A \textbf{graph} is a pair \(\GG = (V,E)\) of a finite set of vertices \(V\) with a relation \(E \subseteq [V]^2\), i.e.\ the 2-element subsets of \(V\) indicating the edges between vertices. The edges are thus unordered pairs \(\{v,w\}\), i.e.\ \(\{v,w\} = \{w,v\}\). This defines \(\GG\) as a simple graph without loops. To denote that \(v\) is a vertex and \(\{v,w\}\) an edge of \(\GG\) we simply write \(v \in \GG\) and \( \in \GG\), respectively.
\end{definition}

Two vertices \(v\) and \(w\) of a graph \(\GG\) are \textbf{adjacent} if \(\{v,w\} \in \GG\). The \textbf{degree} of a vertex \(v\) is given by the cardinality \(|\{w \ | \ \{v,w\} \in \GG\}|\). If all the vertices of \(\GG\) have the same degree \(k\), then \(\GG\) is \textbf{\(k\)-regular}. A \textbf{\(k\)-clique} in a graph is a collection of \(k\) vertices whose induced subgraph is a complete graph, i.e.\ all the vertices in the clique are adjacent.

In the following definition we use the same symbol for a directed graph as in the previous definition for a graph. In the rest of the paper we will always make clear whether we are referring to a graph or a directed graph.

\begin{definition}\label{def:digraph}
	A \textbf{directed graph} (digraph) is a pair \(\GG = (V,E)\) of a finite set of vertices \(V\) and relation \(E \subseteq (V \times V)/\Delta_V\), where \(\Delta_V = \{(v,v) \, | \, v \in V\}\). The relation \(E\) is the set of directed edges between vertices. The edges are unique ordered pairs \((v,w)\), but we allow \textbf{reciprocal edges} \((v,w)\) and \((w,v)\) in \(\GG\). This defines \(\GG\) as a simple directed graph without loops. To denote that \(v\) is a vertex and \((v,w)\) an edge of \(\GG\) we simply write \(v \in \GG\) and \((v,w) \in \GG\), respectively.
\end{definition}

\begin{definition}\label{def:abstract_simplicial complex}
	An \textbf{abstract simplicial complex} on a vertex set \(V\) is a collection \(K\) of non-empty finite subsets \(\sigma \subseteq V\) that is closed under taking non-empty subsets: if \(\sigma \in X\) and \(\tau \subseteq \sigma\) is non-empty then \(\tau \in X\). The subsets are called simplices of \(K\).  
\end{definition}

From now on we drop the word abstract and just talk about simplicial complexes. The following list records notations related to simplices and complexes that we use in this paper.
\begin{center}
	\renewcommand{\arraystretch}{1.5}
	\begin{tabular}{ c | p{8.5cm} }
		Notation & Definition\\ \hline
		\(\sigma \in K\) & \(\sigma\) is a simplex in a complex \(K\). \\ \hline
		\(K_q\) & the set of simplices of \(K\) with dimension  \(\geq q\). \\ \hline
		\(\text{Vert}(K)\), \(\text{Vert}(\sigma)\) & The sets of vertices of \(K\) and \(\sigma\), respectively. \\ \hline
		\(\text{dim}(\sigma)\) & \(|\text{Vert}(\sigma)| - 1\), dimension of \(\sigma\). If equal to \(k\), then \(\sigma\) is a \(k\)-simplex. \\ \hline
		\(\text{dim}(K)\) & The dimension of \(K\) = the dimension of its highest-dimensional simplex. \\ \hline
		\(\tau \hookrightarrow \sigma\) & \(\tau\) is a face of \(\sigma\), i.e.\ \(\tau \subseteq \sigma\). We use the convention that every simplex is a face of itself. Proper face \(\tau\) has dimension strictly less than that of \(\sigma\).		
	\end{tabular}
\end{center}

A simplex is \textbf{maximal} with respect to inclusion if it is not a face of another simplex. The \textbf{\(k\)-skeleton} of a simplicial complex \(K\) is the subcomplex induced by all simplices \(s \in K\) with \(\dim(s) \leq k\). 

\textbf{Flag or clique complex} is a standard way of constructing a simplicial complex from a graph: the \(k\)-simplices are the \((k+1)\)-cliques in the graph. As opposed to graphs and flag complexes, directed graphs and directed analogs of their simplicial complexes are less studied. Directed flag complex is a natural generalisation and the construction we are mostly interested in. By an ordered set we mean an ordered tuple with all entries distinct.

\begin{definition}\label{def:ordered_simplicial complex}
	An \textbf{abstract ordered simplicial complex} on a vertex set \(V\) is a collection of non-empty finite ordered subsets \(\sigma \subseteq V\) that is closed under taking non-empty ordered subsets. All the notations in the above table apply to ordered simplicial complexes.
\end{definition}

\begin{definition}\label{def:directed_flag_complex}
	Let \(\GG = (V,E)\) be a directed graph. The \textbf{directed flag complex} \(\dfl{\GG}\) is the ordered simplicial complex whose \(k\)-simplices are all totally ordered \((k+1)\)-cliques, i.e. \((k+1)\)-tuples \(\sigma = (v_0,v_1,\dots,v_k)\), such that \(v_i \in V\) for all \(i\), and \((v_i,v_j) \in E\) for all \(i < j\). 
\end{definition}
Note that any ordered pair \((v_0,v_i)\) in \(\sigma\), for \(v_0\) fixed, is a directed edge in \(\GG\), and similarly for any ordered pair \((v_i,v_k)\) for \(v_k\) fixed. Therefore the vertex \(v_0\) is called the \textbf{source} of \(\sigma\) and the vertex \(v_k\) is called the \textbf{sink} of \(\sigma\), and simplices obtain a coherent directionality \(v_0 \ra v_k\).

Strictly, ordered simplicial complexes and directed flag complex are examples of an abstract notion of a \textbf{semi-simplicial set} (also known as a \(\Delta\)-set). Our construction of the directed extension of \(Q\)-analysis relies on face maps of semi-simplicial sets so we will recall the definition, see \cite{Friedman2012} for a more in-depth discussion.
\begin{definition}\label{def:semi_simplicial_set}
	A semi-simplicial set \(X\) consists of
	\begin{enumerate}
		\item a sequence of sets \(X_0,X_1,\dots\), and
		\item for each \(n \geq 0\) and for each \(0 \leq i \leq n+1\) \textbf{face maps} \(d_i \colon X_{n+1} \ra X_n\) such that \(d_i d_j = d_{j-1} d_i\) whenever \(i < j\).
	\end{enumerate}
\end{definition}

The sets \(X_n\) are just abstract sets connected by maps \(d_i\). We will not need this abstract picture and will continue to take the sets \(X_n\) as  sets of \(n\)-simplices of an ordered simplicial complex, where the face maps contain the information on how the simplices are attached to each other. 

Note that a simplicial complex has the property that the intersection \(\sigma \cap \tau\) of two simplices is either empty or is a common face of both \(\sigma\) and \(\tau\). Semi-simplicial sets are more general than simplicial complexes since the simplices are not given solely in terms of their vertices. This makes a distinctive difference. Indeed, consider the directed graph below whose 1-dimensional directed flag complex is the graph itself.

\begin{center}
	\begin{tikzpicture}
	\tikzstyle{point}=[circle,thick,draw=black,fill=black,inner sep=0pt,minimum width=2pt,minimum height=2pt]
	\tikzstyle{arc}=[shorten >= 4pt,shorten <= 4pt,->, thick]
	
	\node [point] at (0,0) {};
	\node [point] at (3,0) {};
	\node[left] at (0,0) {0};
	\node[right] at (3,0) {1};
	
	\draw [arc] (0,0) to [bend left=40] (3,0);
	\draw [arc] (3,0) to [bend left=40] (0,0);	
	\end{tikzpicture}
\end{center}

\noindent The vertex set \(\{0,1\}\) spans two different ordered simplices, \((0,1)\) and \((1,0)\). Note that the ordering need not come from any underlying order of the vertex set, \(0 < 1\) is nonetheless a total order in the first case, \(1 < 0\) in the second. Also note that the intersection of the 1-simplices is \(\{0,1\}\) which is not their common face. We will not elaborate on this further and call directed flag complexes simplicial complexes, bearing in mind the general notion of a semi-simplicial set.

\section[Classical Q-analysis]{Classical \(Q\)-analysis}\label{sec:q_analysis}
The \(Q\)-analytical approach to the structure of simplicial complexes does not seem to be widely known in combinatorics, applied topology or network science. We therefore give in this section a comprehensive and self-contained overview of classical \(Q\)-analysis in the case of unordered simplicial complexes. We follow the early works \cite{Atkin1972,Johnson1981}. A modern textbook containing introduction to \(Q\)-analysis is \cite{Johnson_book}.

\subsection[q-connectivity]{\(q\)-connectivity}\label{sec:q_connectivity}

\begin{definition}
	Two simplices \(\sigma\) and \(\tau\) in a simplicial complex \(K\) are \textbf{\(q\)-near}, if they share a \(q\)-face. 
\end{definition}

\begin{definition}\leavevmode\label{def:q_connectivity}
	\begin{enumerate}
		\item Two simplices \(\sigma\) and \(\tau\) of \(K\) are \textbf{\(q\)-connected}, if there is a sequence of simplices in \(K\),
		\[\sigma=\alpha_0,\alpha_1,\alpha_2,\dots,\alpha_n,\alpha_{n+1}=\tau,\]
		such that any two consecutive ones are \(q\)-near. The sequence of simplices is called a \textbf{\(q\)-connection} between \(\sigma\) and \(\tau\). The \textbf{length} of connection is \(n+1\). 
		\item The complex \(K\) is \textbf{\(q\)-connected} if any two simplices in \(K\) of dimension greater than or equal to \(q\) are \(q\)-connected. 
	\end{enumerate}
\end{definition}
\begin{remark}
	Note that the notion of being \(q\)-connected is different from the notion of a topological space \(X\) being \(n\)-connected when all the homotopy groups \(\pi_k(X)\) vanish for \(k \le n\). Indeed, it is easy to come up with an example of a 1-connected simplicial complex, in the language of this paper, with non-vanishing fundamental group. 
\end{remark}
From now on we will write the sequence of a \(q\)-connection as \((\sigma\alpha_1\alpha_2\dots\alpha_n\tau)\). Obviously for a \(q\)-connection to exist between \(\sigma\) and \(\tau\) they have to be of dimension greater than or equal to \(q\). A simplicial complex \(K\) can only have connections up to \(\text{dim}(K)\). The fundamental properties of \(q\)-connectivity are collected in the next proposition.

\begin{proposition}\leavevmode\label{prop:q_connectivity_basic_properties}
	\begin{enumerate}
		\item Every \(q\)-simplex is \(q\)-connected to itself with a connection of length 1.
		\item If \(\sigma\) and \(\tau\) are \(q\)-connected, then they are \(p\)-connected for any \(p < q\).
		\item If \(\sigma\) is maximal with respect to inclusion and \(\text{dim}(\sigma) = q\), then \(\sigma\) is \(q\)-connected only to itself.
		\item If \(\sigma\) and \(\tau\) are \(q\)-connected, any of their \(p\)-faces are \(p\)-connected for \(p < q\).
	\end{enumerate}
\end{proposition}
\begin{proof}\leavevmode
	\begin{enumerate}
		\item Since by definition every simplex is a face of itself, then the trivial sequence \((\sigma\sigma)\) is a \(q\)-connection of length 1.
		\item Since every consecutive pair of simplices in the \(q\)-connection between \(\sigma\) and \(\tau\) shares a common \(q\)-face \(\alpha\), they also share a \(p\)-face of \(\alpha\) for \(p < q\).
		\item Let \(\sigma\) be maximal with respect to inclusion and \(\text{dim}(\sigma) = q\). For \(\sigma\) to be \(q\)-connected to another simplex it has to share a \(q\)-face with another simplex. Since \(\sigma\) is not included in any other simplex and is of dimension \(q\), it can only share a \(q\)-face with itself.
		\item Let \((\sigma\alpha_1\alpha_2\dots\alpha_n\tau)\) be the \(q\)-connection between \(\sigma\) and \(\tau\). Adjoin \(p\)-faces \(f_\sigma\) and \(f_\tau\) of \(\sigma\) and \(\tau\), respectively, to get a sequence \((f_\sigma \sigma\alpha_1\alpha_2\dots\alpha_n\tau f_\tau)\). The \(p\)-faces are \(p\)-near to their respective simplices and by statement 2. of the proposition, the rest of the sequence provides a \(p\)-connection for any \(p < q\). 
	\end{enumerate}
\end{proof}

Since the highest dimensional simplices of \(K\) are maximal, they are \(\text{dim}(K)\)-connected only to themselves by the third point in Proposition \ref{prop:q_connectivity_basic_properties}. The next statement forms the foundation of \(Q\)-analysis.

\begin{theorem}\label{prop:q_connectivity_equivalence}
	Let \(K_q\) denote the set of simplices of \(K\) with dimension greater than or equal to \(q\). For any \(\sigma\) and \(\tau\) in \(K_q\), the relation \(\sim_q\) defined by
	\[\sigma \sim_q \tau \ \text{if and only if} \ \sigma \ \text{and} \ \tau \ \text{are} \ q\text{-connected},\]
	is an equivalence relation. The equivalence classes, i.e.\ the elements, of the quotient \(K_q/\!\! \sim_q\) are called the \(q\)-\textbf{connected components} of \(K\).
\end{theorem}
\begin{proof}
	Reflexivity follows by definition since every simplex of dimension \(q\) or higher is \(q\)-connected to itself. For symmetry, if \((\sigma\alpha_1\alpha_2\dots\alpha_n\tau)\) is a \(q\)-connection from \(\sigma\) to \(\tau\), then the reverse sequence obviously is a \(q\)-connection from \(\tau\) to \(\sigma\). For transitivity, if \((\sigma\alpha_1\alpha_2\dots\alpha_n\tau)\) and \((\tau\beta_1\beta_2\dots\beta_n\kappa)\) are \(q\)-connections between \(\sigma\) and \(\tau\), and \(\tau\) and \(\kappa\), respectively, then \((\sigma\alpha_1\alpha_2\dots\alpha_n\tau\beta_1\beta_2\dots\beta_n\kappa)\) is a \(q\)-connection between \(\sigma\) and \(\kappa\).
\end{proof}

\begin{remark}\label{rem:q_category}
	As was noted in the original works of Atkin, the equivalence relation above immediately gives us a small category \(\Gamma_q(K)\) with objects \(K_q\), and unique morphisms \(\sigma \ra \tau\) whenever \(\sigma \sim_q \tau\). Note that \(\Gamma_0(K)\) contains all simplices of \(K\) with all possible connections. Since \(K_{q+1} \subseteq K_q\) and if \(\sigma \sim_{q+1} \tau\) then \(\sigma \sim_q \tau\), we have a sequence of subcategories \(\Gamma_n(K) \subseteq \dots \subseteq \Gamma_1(K) \subseteq \Gamma_0(K)\), where \(n = \text{dim}(K)\).
\end{remark}

We recall the notion of a pseudomanifold, which can be formulated in terms of \(q\)-connectivity \cite{Spanier}.

\begin{definition}\label{def:pseudomanifold}
	A simplicial complex \(K\) is an \(n\)\textbf{-pseudomanifold} if all the maximal simplices are \(n\)-simplices, each \((n-1)\)-simplex is a face of exactly two \(n\)-simplices, and any two \(n\)-simplices are \((n-1)\)-connected.
\end{definition}

\begin{definition}
	A simplicial complex \(K\) is an \(n\)\textbf{-pseudomanifold with boundary} if all the maximal simplices are \(n\)-simplices, each \((n-1)\)-simplex is a face of at most two \(n\)-simplices, and any two \(n\)-simplices are \((n-1)\)-connected.
\end{definition}
The boundary of an \(n\)-pseudomanifold \(K\), denoted \(\partial K\), is the subcomplex of \(K\) given by those \((n-1)\)-simplices, each of which is a face of exactly one \(n\)-simplex.


\subsection[Q-vectors of simplicial complexes]{\(Q\)-vectors of simplicial complexes}
The \textbf{\(Q\)-analysis} of a simplicial complex means to find its \(q\)-connectivity classes, i.e.\ the \(q\)-connected components for \(0 \leq q \leq \text{dim}(K)\). This information is summarized with various structure vectors. Let \(Q_q\) denote the number of \(q\)-connectivity classes and let \(n = \text{dim}(K)\).
\begin{definition}
	The \textbf{first structure vector} of \(K\) is the tuple
	\[\mathbf{Q}(K) = (Q_n, Q_{n-1}, \dots, Q_0).\]
\end{definition} 
Since 0-connectivity corresponds to path connectedness, \(Q_0\) is equal to the Betti number \(\beta_0\), the number of connected components of \(K\). The elements of \(\mathbf{Q}(K)\) are therefore generalizations of \(\beta_0\) to higher-dimensional connectivity. As noted after Proposition \ref{prop:q_connectivity_basic_properties}, \(Q_n\) is the number of highest-dimensional simplices of \(K\). All the elements of \(\mathbf{Q}(K)\) are \(\geq 1\) since every simplex is at least connected to itself (Proposition \ref{prop:q_connectivity_basic_properties}).

Other structure vectors have also been defined in the literature. The \textbf{second structure vector} of \(K\) is defined as 
\[\mathbf{N}(K) = (|K_n|, |K_{n-1}|, \dots, |K_0|),\]
where \(|K_n|\) is the number of simplices of dimension \(n\) or higher. The \textbf{third structure vector} \cite{TrafficJamming}, or the reduced structure vector, is then defined by
\[\mathbf{T}(K) = (1-Q_n/|K_n|, 1-Q_{n-1}/|K_{n-1}|, \dots, 1-Q_0/|K_0|).\]

\begin{remark}
	Recall that a map \(\phi \colon K \ra L\) between simplicial complexes is a simplicial isomorphism if it is bijective as a function \(\phi \colon \text{Vert}(K) \ra \text{Vert}(L)\), and whenever \(\{v_0,\dots,v_n\}\) is a simplex in \(K\) then \(\{\phi(v_0),\dots,\phi(v_n)\}\) is a simplex in \(L\). It is then plain that simplicial isomorphisms preserve \(q\)-connectivities between simplices and hence the structure vectors, and \(Q\)-analysis is a simplicial isomorphism invariant.
\end{remark}

Let \(\mathbf{1}_n\) denote a vector of length \(n+1\) with all components equal to 1.

\begin{proposition}\label{prop:simplex_vector_characterization}
	\(K\) is an \(n\)-dimensional simplex if and only if \(\mathbf{Q}(K) = \mathbf{1}_n\).
\end{proposition}
\begin{proof}
	\(\Rightarrow \colon\) Since \(K\) is a simplex, the 0th and \(n\)th elements of \(\mathbf{Q}(K)\) are 1. Any proper face of dimension \(d < n\) is \(d\)-connected to any other face of dimension \(d\) or higher through \(K\). There is hence only 1 \(d\)-connected component for \(n > d > 0\) and the claim follows.
	
	\noindent \(\Leftarrow \colon\) Assume \(\mathbf{Q}(K) = \mathbf{1}_n\) and \(K\) is an \(n\)-dimensional complex which is not an \(n\)-simplex. We can assume \(K\) is connected because \(Q_0 = 1\). Therefore \(K\) contains a maximal simplex of some dimension \(0 < q \leq n\). By Proposition \ref{prop:q_connectivity_basic_properties} this simplex is \(q\)-connected only to itself, \(Q_q > 1\) and we arrive at a contradiction.
\end{proof}

\(Q\)-analysis associates to \(K\) one more vector, whose definition should be apparent after Proposition \ref{prop:simplex_vector_characterization}, which essentially says that simplices are fully \(q\)-connected for any \(q\). 
\begin{definition}
	The \textbf{obstruction vector} of \(K\), for \(\text{dim}(K) = n\), is defined by
	\[\mathbf{\hat{Q}}(K) = \mathbf{Q}(K) - \mathbf{1}_n.\]
\end{definition}
The obstruction vector is then a measure of how much \(K\) deviates from being an \(n\)-simplex. For example, for a contractible complex which is not a simplex, all Betti numbers beyond \(\beta_0\) are zero. The vectors \(\mathbf{Q}(K)\) and \(\mathbf{\hat{Q}}(K)\) are therefore more sensitive to the actual combinatorial structure than simplicial homology can detect.

We can use \(q\)-analysis to measure how special a simplex is within a complex. The \textbf{eccentricity} of a simplex \(\sigma\) is defined as
\[\text{ecc}(\sigma) = \frac{\text{dim}(\sigma) - \check{q}}{\check{q} + 1},\]
where \(\check{q}\) is the greatest value of \(q\) for which \(\sigma\) is \(q\)-connected to another simplex which is not a face of \(\sigma\). Note that \(\check{q} + 1\) is the greatest number of vertices \(\sigma\) uses to share in connection to another simplex, while \(\text{dim}(\sigma) - \check{q}\) is the number of vertices left independent from any connection. Eccentricity therefore accords with the intuition of how isolated a simplex is from \(q\)-connectivity point of view. For a \(\sigma\) that is \(\text{dim}(\sigma)\)-connected to another simplex, \(\text{ecc}(\sigma) = 0\).

\begin{example}	
	The figure below depicts two 3-dimensional complexes \(K\) and \(K'\). Their first and second structure vectors, as well as the eccentricities of one of their 3-simplices are also shown. The change in the \(\mathbf{Q}\)-vector when going from \(K\) to \(K'\) indicates shift to higher connectivity by filling the combinatorial 1-dimensional hole in \(K\). This can also be measured by the obstruction vector \(\mathbf{\hat{Q}}(K') = (1,4,0,0)\).
	\begin{center}
		\begin{tikzpicture}
		\tikzstyle{point}=[circle,thick,draw=black,fill=black,inner sep=0pt,minimum width=2pt,minimum height=2pt]
		
		\draw (0,0) -- (1.5,0) -- (0.75,-1) --cycle; 
		
		\draw[fill=magenta,opacity=0.5] (0,0) -- (-0.75,-1) -- (0.75,-1) --cycle; 
		
		\draw[fill=magenta,opacity=0.5] (0,0) -- (0.75,0.4) -- (1.5,0) --cycle;
		\draw[fill=magenta,opacity=0.5] (0,0) -- (0.75,1) -- (0.75,0.4) --cycle;
		\draw[fill=magenta,opacity=0.5] (0.75,1) -- (0.75,0.4) -- (1.5,0) --cycle;
		\draw (0.75,1) -- (0.75,0.4);			
		\draw[fill=magenta,opacity=0.5] (0,0) -- (0.75,1) -- (1.5,0) --cycle;
		
		\draw[fill=magenta,opacity=0.5] (1.5,0) -- (0.75,-1) -- (2.25,-1) --cycle;
		
		\draw[fill=magenta,opacity=0.4] (0.75,-1) -- (2.25,-1) -- (1.5,-1.4) --cycle;
		\draw[fill=magenta,opacity=0.4] (0.75,-1) -- (1.5,-1.4) -- (1.5,-2) --cycle;
		\draw[fill=magenta,opacity=0.4] (2.25,-1) -- (1.5,-1.4) -- (1.5,-2) --cycle;
		\draw (0.75,1) -- (0.75,0.4);		
		\draw[fill=magenta,opacity=0.4] (0.75,-1) -- (2.25,-1) -- (1.5,-2) --cycle;
		
		\node[point] at (0,0) {};
		\node[point] at (-0.75,-1) {};
		\node[point] at (0.75,-1) {};
		\node[point] at (0.75,0.4) {};
		\node[point] at (1.5,0) {};
		\node[point] at (0.75,1) {};
		\node[point] at (2.25,-1) {};
		\node[point] at (1.5,-1.4) {};
		\node[point] at (1.5,-2) {};
		
		\node at (-0.5,0.5) {\(\text{ecc} = 3\)};
		\node[right] at (-0.5,-2.5) {\(\mathbf{Q}(K) = (2,4,3,1)\)};
		\node[right] at (-0.5,-3) {\(\mathbf{N}(K) = (2,12,29,38)\)};
		
		\draw[fill=magenta,opacity=0.5] (6,0) -- (7.5,0) -- (6.75,-1) --cycle; 
		
		\draw[fill=magenta,opacity=0.5] (6,0) -- (5.25,-1) -- (6.75,-1) --cycle; 
		
		\draw[fill=magenta,opacity=0.5] (6,0) -- (6.75,0.4) -- (7.5,0) --cycle;
		\draw[fill=magenta,opacity=0.5] (6,0) -- (6.75,1) -- (6.75,0.4) --cycle;
		\draw[fill=magenta,opacity=0.5] (6.75,1) -- (6.75,0.4) -- (7.5,0) --cycle;
		\draw (6.75,1) -- (6.75,0.4);			
		\draw[fill=magenta,opacity=0.5] (6,0) -- (6.75,1) -- (7.5,0) --cycle;
		
		\draw[fill=magenta,opacity=0.5] (7.5,0) -- (6.75,-1) -- (8.25,-1) --cycle;
		
		\draw[fill=magenta,opacity=0.4] (6.75,-1) -- (8.25,-1) -- (7.5,-1.4) --cycle;
		\draw[fill=magenta,opacity=0.4] (6.75,-1) -- (7.5,-1.4) -- (7.5,-2) --cycle;
		\draw[fill=magenta,opacity=0.4] (8.25,-1) -- (7.5,-1.4) -- (7.5,-2) --cycle;
		\draw (6.75,1) -- (6.75,0.4);		
		\draw[fill=magenta,opacity=0.4] (6.75,-1) -- (8.25,-1) -- (7.5,-2) --cycle;
		
		\node[point] at (6,0) {};
		\node[point] at (5.25,-1) {};
		\node[point] at (6.75,-1) {};
		\node[point] at (6.75,0.4) {};
		\node[point] at (7.5,0) {};
		\node[point] at (6.75,1) {};
		\node[point] at (8.25,-1) {};
		\node[point] at (7.5,-1.4) {};
		\node[point] at (7.5,-2) {};
		
		\node at (5.5,0.5) {\(\text{ecc} = 1\)};
		\node[right] at (5.5,-2.5) {\(\mathbf{Q}(K') = (2,5,1,1)\)};
		\node[right] at (5.5,-3) {\(\mathbf{N}(K') = (2,13,30,39)\)};	
		\end{tikzpicture}
	\end{center}
\end{example}

As in Proposition \ref{prop:simplex_vector_characterization} for an \(n\)-simplex, we know the first structure vectors of pseudomanifolds.
\begin{proposition}\label{prop:pseudomanifold_vector_characterization}
	Let \(K\) be an \(n\)-pseudomanifold, with or without boundary, and let its number of \(n\)-simplices be \(t\). Then \(\mathbf{Q}(K) = (t,1,1,\dots,1)\).
\end{proposition}
\begin{proof}
	By definition any \(k\)-simplex for \(k < n\) is a face of some \(n\)-simplex. By definition as well, all the \(n\)-simplices are \((n-1)\)-connected. It then follows from statements 2. and 4. in Proposition \ref{prop:q_connectivity_basic_properties} that \(Q_k = 1\) for any \(k < n\). Each of the \(t\) maximal \(n\)-simplices is an \(n\)-connected component in itself so \(Q_n = t\).
\end{proof}

The illustration below is a contradiction to the necessity of Proposition \ref{prop:pseudomanifold_vector_characterization}. Both complexes have first structure vectors \((3,1,1)\) and 3 maximal 2-simplices but the right one is not a pseudomanifold since the middle 1-simplex is a face of more than two 2-simplices.
\begin{center}
	\begin{tikzpicture}
	\tikzstyle{point}=[circle,thick,draw=black,fill=black,inner sep=0pt,minimum width=2pt,minimum height=2pt]
	\draw[fill=magenta,opacity=0.5] (0,0) -- (1.5,0) -- (0.75,-1) --cycle; 
	
	\draw[fill=magenta,opacity=0.5] (0,0) -- (-0.75,-1) -- (0.75,-1) --cycle; 
	
	\draw[fill=magenta,opacity=0.5] (1.5,0) -- (0.75,-1) -- (2.25,-1) --cycle;
	
	\node[point] at (0,0) {};
	\node[point] at (-0.75,-1) {};
	\node[point] at (0.75,-1) {};
	\node[point] at (1.5,0) {};
	\node[point] at (2.25,-1) {};
	\draw[fill=magenta,opacity=0.5] (5,0) -- (5,-1) -- (4.2,0.4) --cycle;
	
	\draw[fill=magenta,opacity=0.5] (5,0) -- (5,-1) -- (4,-1.5) --cycle; 
	
	\draw[fill=magenta,opacity=0.5] (5,0) -- (5,-1) -- (6,-1.5) --cycle;
	
	\node[point] at (5,0) {};
	\node[point] at (5,-1) {};
	\node[point] at (4.2,0.4) {};
	\node[point] at (4,-1.5) {};
	\node[point] at (6,-1.5) {};
	
	\end{tikzpicture}
\end{center}

\subsection[Algorithm for incidence structures]{\(Q\)-analysis of incidence structures}
One of the initial motivations to define the concept of \(Q\)-analysis came from a geometric study of relational structures, particularly in applications to social contexts. Let \(X\) and \(Y\) be finite sets. Recall that a relation \(R\) between \(X\) and \(Y\) is a subset of the product \(X \times Y\). When \((x,y) \in R\) we write \(x R y\). The \textbf{incidence matrix} of a relation is the \(|X| \times |Y|\) binary matrix \(\Lambda(R)\), where \(\Lambda(R)_{ij} = 1\) if and only if \(x_i R y_j\). The inverse relation \(R^{-1}\) is a subset of \(Y \times X\) such that \(y R^{-1} x\) if and only if \(x R y\). The incidence matrix of the inverse relation is \(\Lambda(R^{-1}) = \Lambda(R)^T\). We also just write \(\Lambda\) for an incidence matrix.

To a relation \(R\) we can associate two simplicial complexes. The complex \(K_X(Y,R)\) has as its vertices the set \(Y\) and as simplices all subsets \(\sigma \subset Y\) whenever there is an \(x \in X\) such that \(x R y\) for all \(y \in \sigma\). Simplices of \(K_X(Y,R)\) can be read off from the rows of \(\Lambda(R)\) as those subsets of elements with value 1. Likewise, the complex \(K_Y(X,R)\) has as its vertices the set \(X\) and as simplices all subsets \(\sigma \subset X\) whenever there is a \(y \in Y\) such that \(y R^{-1} x\) for all \(x \in \sigma\). The columns of \(\Lambda(R)\) indicate the simplices in \(K_X(Y,R)\) as those subsets of elements with value 1. The classic work of Dowker \cite{Dowker1952} established that the homology groups of \(K_X(Y,R)\) and \(K_Y(X,R)\) are isomorphic. The \(Q\)-analyses of these complexes, however, are different as noted by Atkin \cite{Atkin1972}.

The \(Q\)-analytic information of \(R\) can be computed from the product of incidence matrices \(\Lambda\Lambda^T\). The algorithm is evident by seeing how the element \(\lambda_{ij}\) of the product is given:

\[
\begin{blockarray}{ccccc}
& y_1 & y_2 & .. & y_n  \\
\begin{block}{c[cccc]}
.. &  &  &  &  \\
x_i & - & - & - & - \\
.. &  &  &  &  \\
x_m &  &  &  & \\
\end{block}
\end{blockarray}
\begin{blockarray}{ccccc}
.. & x_j  & .. & x_m &\\
\begin{block}{[cccc]c}
& \mid &  &  & y_1 \\
& \mid &  &  & y_2 \\
& \mid &  &  & .. \\
& \mid &  &  & y_n \\
\end{block}
\end{blockarray}
\]
The row \(x_i\) and column \(x_j\) denote simplices in the complex \(K_X(Y,R)\). As binary vectors, their inner product \(\lambda_{ij}\) is the number of vertices the corresponding simplices share, i.e.\ \(\lambda_{ij} = \text{dim}(x_i \cap x_j) + 1.\) The matrix \(\Lambda\Lambda^T - \mathbf{1}_{m,m}\), where \(\mathbf{1}_{m,m}\) denotes the \(m \times m\) matrix of all ones, is therefore the matrix of dimensions of shared faces between simplices in \(K_X(Y,R)\). Since a simplicial complex can be put into an incidence matrix with columns labeled by vertices and rows representing all higher-dimensional simplices, the above algorithm computes the \(q\)-connectivity structure of a complex. 

\begin{definition}\label{def:q_graph}
	The \(q\)-\textbf{graph} of a simplicial complex \(K\) has as its vertices the simplices in \(K_q\) and edges between pairs of \(q\)-near simplices. 
\end{definition}

The \(q\)-graph and \(q\)-connected components of \(K\) can now be computed by the following steps:
\begin{enumerate}
	\item Form the \(m \times n\) incidence matrix \(\Lambda\) whose rows are labeled by the \(m\) simplices in \(K_q\) and columns are labeled by the \(n+1\) vertices of the elements in \(K_q\).
	\item Compute the matrix \(\Lambda\Lambda^T - \mathbf{1}_{m,m}\) and set to one all \(\lambda_{ij} \geq q\) and to zero all other entries.
	\item The obtained matrix represents the \(q\)-nearness relation and is the adjacency matrix of the \(q\)-graph. The connected components of the graph are the equivalence classes of the \(q\)-connectivity relation, which can be obtained by standard graph algorithms.
\end{enumerate} 
Note that the above algorithm is not efficient for a large number of simplices since it requires initialising a large sparse matrix. It is then more efficient to construct the \(q\)-graph edge by edge through pairwise comparison of simplices. See \cite{LiKwong2009} for another algorithm for extracting the first structure vector and \(q\)-connected components from the product \(\Lambda\Lambda^T\). Note that the matrix is symmetric so it suffices to only look at the upper/lower triangle to extract the \(q\)-analytic information. The following statement gives characterisations of pseudomanifolds in terms of \(q\)-graphs.

\begin{proposition}\leavevmode\label{prop:pseudomanifolds_q_graphs}
	\begin{enumerate}
		\item The transitive closure of the \((n-1)\)-graph of an \(n\)-pseudomanifold, with or without boundary, is a complete graph.
		\item Let \(\GG\) be the \((n-1)\)-graph of an \(n\)-pseudomanifold. Then the subgraph induced by \(n\)-simplices is \((n+1)\)-regular.
	\end{enumerate}
	
\end{proposition}
\begin{proof}\leavevmode
	\begin{enumerate}
		\item The vertices are all the \((n-1)\)- and \(n\)-simplices. By definition of pseudomanifolds there is a path between any two \(n\)-simplices. Similarly any \((n-1)\)-simplex is near to some \(n\)-simplex and by the above path there is a path to any \((n-1)\)- or \(n\)-simplex. Taking the transitive closure therefore gives a complete graph.
		\item By Definition \ref{def:pseudomanifold} any \((n-1)\)-simplex as a vertex of \(\GG\) has degree 2, being adjacent to exactly two \(n\)-simplices \(x\) and \(y\). This implies that \(x\) and \(y\) are \((n-1)\)-near and there is an edge \(\{x,y\}\). Hence, every \((n-1)\)-face of an \(n\)-simplex \(x\) corresponds to an edge \(\{x,y\}\) to some \(n\)-simplex \(y\). Moreover, every \(y\) is different for a different \((n-1)\)-face since otherwise \(K\) would not be a pseudomanifold by the above degree condition. As the number of \((n-1)\)-faces of an \(n\)-simplex is \(n+1\), the subgraph induced by the \(n\)-simplices is \((n+1)\)-regular. 
	\end{enumerate}
	
\end{proof}
The above proposition gives ways to computationally recognise pseudomanifolds, compare also to Proposition \ref{prop:pseudomanifold_vector_characterization}. Of course since \(q\)-connectivity is an equivalence relation the transitive closures of \(q\)-graphs are complete. But for a general simplicial complex there is no guarantee that the graph stays connected. 

There are two well-known notions, one from graph theory and another from (topological) combinatorics: \(k\)-clique community \cite{clique_communities} and face poset \cite{Barmak_book}. We record here the observation, not seeming to appear in the literature, that these two notions arise as subrelations of \(q\)-connectivity.
\begin{definition}
	Let \(\GG\) be a graph and \(k \geq 2\). Two \(k\)-cliques in \(\GG\) are connected if there is a sequence of \(k\)-cliques of \(\GG\) such that any two consecutive cliques share \(k-1\) vertices. A \textbf{\(k\)-clique community} of \(\GG\) is a maximal set of pairwise connected \(k\)-cliques.
\end{definition}

\begin{proposition}
	A \(k\)-clique community of a graph \(\GG\) is a \((k-2)\)-connected component in the \((k-1)\)-skeleton of the flag complex of \(\GG\).
\end{proposition}
\begin{proof}
	The \(k\)-cliques of \(\GG\) are the \((k-1)\)-simplices in the flag complex of \(\GG\), and therefore are contained in the \((k-1)\)-skeleton. For \((k-2)\)-connectivity of the skeleton we look at the set of simplices \(K_{k-2}\), which now consists only of \((k-1)\)- and \((k-2)\)-simplices. Two \(k\)-cliques being connected in \(\GG\) means exactly that the corresponding \((k-1)\)-simplices in the \((k-1)\)-skeleton are \((k-2)\)-connected. Taking the vertices of the \((k-1)\)-simplices in the connected components back to cliques of \(\GG\) recovers the \(k\)-clique communities.
\end{proof}

\begin{definition}
	Let \(K\) be a simplicial complex. The \textbf{face poset} of \(K\) has simplices of \(K\) as elements and inclusion of simplices as the partial order relation.
\end{definition}
In the following proposition we specifically consider the face poset as a relation \(F=\{(\sigma,\tau) \ | \ \sigma \hookrightarrow \tau\}\). Recall the union of relations. Let \(R\) and \(S\) be relations on sets \(X\) and \(Y\), respectively. Then the union \(R \cup S\) on the set \(X \cup Y\) is defined to be \(\{(x,x') \ | \ xRx'\} \cup \{(y,y') \ | \ ySy'\}\).
\begin{proposition}\label{prop:face_poset_inside_q_connection}
	Let \(K\) be a simplicial complex and \(F\) its face poset. Let \(R_q\) denote the relation of being \(q\)-near on the set of simplices \(K_q\). Then \(F \subset \bigcup_{q=0}^{\dim(K)} R_q\).
\end{proposition}
\begin{proof}
	We can construct the Hasse diagram of \(H(F)\) from \(\bigcup_{q=0}^{\dim(K)} R_q\) which determines \(F\) completely. For a simplex \(\sigma\), all pairs \((\sigma,\tau)\) are in \(H(F)\) where \(\sigma \hookrightarrow \tau\) and there is no \(\alpha\) such that \(\sigma \hookrightarrow \alpha \hookrightarrow \tau\). This means that \(\sigma\) and \(\tau\) are \(\dim(\sigma)\)-near and the pairs \((\sigma,\tau)\) are those in \(R_{\dim(\sigma)}\) where \(\dim(\tau) = \dim(\sigma) + 1\).
\end{proof}

The preceding propositions show that if one computes the \(Q\)-analytic information, in particular the \(q\)-graphs, then one also obtains the \(k\)-clique communities and the face poset by appropriate restrictions. It is well known that the face poset has the same homotopy type as the simplicial complex it is associated to \cite{Barmak_book}. The relational structure of \(q\)-connection therefore contains the homotopical information of the complex \(K\). 

\section[Directed Q-analysis]{Directed \(Q\)-analysis}\label{sec:directed_q_analysis}
The main construction of this paper is now to extend \(q\)-connectivity to ordered simplicial complexes such as directed flag complexes built out of directed graphs. Note that the notion of connectivity in the previous section is independent of directionality of the simplices: any \(q\)-face shared between directed simplices is again a directed simplex and as such we can study \(q\)-connectivity of directed complexes. However, our motivation for the construction in this section is the observation that standard \(Q\)-analysis fails to distinguish connectivity features occurring in directed flag complexes.
\begin{example}\label{ex:q_connectivity_problem}
	The directed graphs, or equally 1-dimensional directed flag complexes, depicted below have the same \(q\)-connectivity structure as undirected graphs/simplicial complexes. The maximal 1-simplices each form their own 1-connected components and both complexes are connected. The first structure vectors in both cases are therefore (3,1).
	\begin{center}
		\begin{tikzpicture}
		\tikzstyle{point}=[circle,thick,draw=black,fill=black,inner sep=0pt,minimum width=2pt,minimum height=2pt]
		\tikzstyle{arc}=[shorten >= 8pt,shorten <= 8pt,->, thick]
		
		\node[] at (-1.2,-1) {1};
		\node[] at (1.2,-1) {2};
		\node[] at (0,1) {0};
		
		\draw[arc] (0,1) to (-1.2,-1);
		\draw[arc] (-1.2,-1) to (1.2,-1);
		\draw[arc] (1.2,-1) to (0,1);
		
		\node[] at (5,-0.3) {a};
		\node[] at (3.8,-1) {b};
		\node[] at (6.2,-1) {c};
		\node[] at (5,1) {d};		
		
		\draw[arc] (3.8,-1) to (5,-0.3);
		\draw[arc] (6.2,-1) to (5,-0.3);
		\draw[arc] (5,1) to (5,-0.3);
		\end{tikzpicture}
	\end{center}
\end{example}

Standard \(q\)-connectivity is only sensitive to shared \(q\)-faces. In the Example \ref{ex:q_connectivity_problem} the structure vectors are therefore unable to capture the very different connectivity arising from directed edges. Our motivation then is to take into account the directionality of simplices. Instead of just sharing a \(q\)-face, we need to impose conditions on \emph{how} the face is shared between two simplices. This is enabled by the face maps (Definition \ref{def:semi_simplicial_set}). If \(\sigma\) is a simplex, then \(d_i(\sigma)\) is a \textbf{face in the \(d_i\)-direction}. Depending on \(i\) and the dimension of \(\sigma\), taking the face in the \(d_i\)-direction might not make sense. We therefore introduce a slightly modified face map.
\begin{definition}\label{def:new_face_map}
	Let \(\sigma\) be an \(n\)-simplex. We denote by \(\widehat{d_i}\) the face map
	\[\widehat{d_i}(\sigma) = 
	\begin{cases}
	(v_0,\dots,\widehat{v_i},\dots, v_n), \ \text{if } i < n,\\
	(v_0,\dots,v_{n-1}, \widehat{v_n}), \ \text{if } i \geq n.
	\end{cases}\]
\end{definition}
The face map \(\widehat{d_i}\) now makes sense in any dimension since it always removes the vertex at position \(\text{min}\{i,\text{dim}(\sigma)\}\).
\begin{remark}
	Note that in the case \(i \geq n\) there is a choice to remove the last vertex. While developing the theory many different definitions were introduced and the above seems to give the simplest and the most natural way to capture the aspects of directed \(q\)-connectivity for our purposes. One may view this as augmenting an \(n\)-simplex with "phantom" vertices, and then the face map \(\widehat{d_i}\), for \(i \geq n\), removes the first actual vertex of \(\sigma\) it can, i.e.\ the last one.
\end{remark}

\begin{definition}\label{def:q_near_directed}
	For an ordered simplicial complex \(K\), let \((\sigma,\tau)\) be an ordered pair of simplices \(\sigma\) and \(\tau\) with \(\text{dim}(\sigma),\text{dim}(\tau) \geq q\). Let \((\widehat{d_i},\widehat{d_j})\) be an ordered pair of face maps. Then \((\sigma,\tau)\) is \textbf{\(q\)-near along \((\widehat{d_i},\widehat{d_j})\)} if either of the following conditions is true:

	\begin{enumerate}
		\item \(\sigma \hookrightarrow \tau,\)
		\item \(\widehat{d_i}(\sigma) \hookleftarrow \alpha \hookrightarrow \widehat{d_j}(\tau),\) for some \(q\)-simplex \(\alpha \in K\).
	\end{enumerate}
\end{definition}
\begin{definition}\label{def:q_connection_directed}
	The ordered pair \((\sigma,\tau)\) of simplices of \(K\) is \textbf{\(q\)-connected along \((\widehat{d_i},\widehat{d_j})\)} if there is a sequence of simplices in \(K\),
	\[\sigma=\alpha_0,\alpha_1,\alpha_2,\dots,\alpha_n,\alpha_{n+1}=\tau,\]
	such that any ordered pair \((\alpha_i,\alpha_{i+1})\) is \(q\)-near along \((\widehat{d_i},\widehat{d_j})\). The sequence of simplices is called a \textbf{\(q\)-connection along \((\widehat{d_i},\widehat{d_j})\)} between \(\sigma\) and \(\tau\). We simply write this connection as \((\sigma \alpha_1 \alpha_2 \dots \alpha_n \tau)\).
\end{definition}
We will call the above connection \((q,\widehat{d_i},\widehat{d_j})\)-connection, when the choices of \(q\) and directions \(\widehat{d_i}\) and \(\widehat{d_j}\) are made, and similarly we say \((q,\widehat{d_i},\widehat{d_j})\)-near. We also simply say \(q\)-near, -connected and -connection when the directed nature along the pair \((\widehat{d_i},\widehat{d_j})\) is clear from the context. 
\begin{proposition}\leavevmode\label{prop:directed_q_connectivity_basic_properties}
	\begin{enumerate}
		\item The pair \((\sigma,\sigma)\) of \(q\)-simplices is \(q\)-connected along any pair \((\widehat{d_i},\widehat{d_j})\).
		\item If \((\sigma,\tau)\) is \(q\)-connected along \((\widehat{d_i},\widehat{d_j})\), then it is \(p\)-connected along \((\widehat{d_i},\widehat{d_j})\) for any \(p < q\).
		\item If \(\sigma\) is maximal with respect to inclusion and \(\text{dim}(\sigma) = q\), then \(\sigma\) is \(q\)-connected only to itself along any pair \((\widehat{d_i},\widehat{d_j})\).
	\end{enumerate}
\end{proposition}
\begin{proof}\leavevmode
	\begin{enumerate}
		\item This follows by definition from \(\sigma \hookrightarrow \sigma\) being \(q\)-near.
		\item Let \((\sigma \alpha_1 \alpha_2 \dots \alpha_n \tau)\) be a \(q\)-connection between \(\sigma\) and \(\tau\). If \(\sigma \hookrightarrow \alpha_1\), then \((\sigma,\alpha_1)\) is automatically \(p\)-near. Otherwise there is a \(q\)-simplex \(\beta\) and a \(p\)-simplex \(\beta'\), \(p < q\), such that \(\widehat{d_i}(\sigma) \hookleftarrow \beta \hookleftarrow \beta' \hookrightarrow \beta \hookrightarrow \widehat{d_j}(\alpha_1).\) As the above happens for any consecutive, ordered pair in the connection \((\sigma \alpha_1 \alpha_2 \dots \alpha_n \tau)\) the assertion follows.
		\item Since \(\sigma\) is maximal it can only be \(q\)-near to itself through \(\sigma \hookrightarrow \sigma\). And since \(\widehat{d_i}(\sigma)\) is of dimension \(q-1\), there cannot be a \(q\)-simplex \(\alpha\) for the relation \(\widehat{d_i}(\sigma) \hookleftarrow \alpha \hookrightarrow \widehat{d_j}(\tau)\) to exist for any \(\tau\).
	\end{enumerate}
\end{proof}
The connection in Definition \ref{def:q_connection_directed} is clearly not symmetric, nor antisymmetric, in general. If \((\sigma\alpha_1\alpha_2\dots\alpha_n\tau)\) and \((\tau\beta_1\beta_2\dots\beta_n\kappa)\) are \(q\)-connections along \((\widehat{d_i},\widehat{d_j})\) between \(\sigma\) and \(\tau\), and \(\tau\) and \(\kappa\), respectively, then we can form a \(q\)-connection between \(\sigma\) and \(\kappa\) by \((\sigma\alpha_1\alpha_2\dots\alpha_n\tau\beta_1\beta_2\dots\beta_n\kappa)\), since \((\tau,\tau)\) is always \(q\)-connected. Hence the directed \(q\)-connection is transitive in the same way as standard \(q\)-connection. By definition we always have reflexivity. The most significant difference from standard \(Q\)-analysis is then encapsulated by the following main result.
\begin{theorem}
	The relation of being \((q,\widehat{d_i},\widehat{d_j})\)-connected is a preorder on \(K_q\).
\end{theorem}
By Definition \ref{def:q_connection_directed} and Proposition \ref{prop:directed_q_connectivity_basic_properties} the \((q,\widehat{d_i},\widehat{d_j})\)-connection is a directed extension of standard \(q\)-connectivity with similar properties. But the directed connection as a preorder relation makes these two approaches quite different. The equivalence relation imposed by standard \(q\)-connectivity associates to a simplicial complex its generalized path components in a canonical way as the connected components of the \(q\)-graph (Definition \ref{def:q_graph}). The preorder structure of \((q,\widehat{d_i},\widehat{d_j})\)-connection associates to an ordered simplicial complex a directed \(q\)-graph, which captures the simplicial connections imposed by the choice of directions \(\widehat{d_i}\) and \(\widehat{d_j}\). We will emphasise this further in Section \ref{sec:dir_struct_sys} where we look at the preorders as finite topological spaces. From the point of view of analysing complexes associated to digraphs this choice offers flexibility to peer into the simplicial structure and we will explore this in Section \ref{sec:applications}.

The analog of the \(q\)-graph in the directed setting will have edges coming only from \((q,\widehat{d_i},\widehat{d_j})\)-nearness of simplices. This is in fact the Hasse diagram form of the connectivity perorders, i.e.\ we do not draw the transitive edges, and we do not draw the reflexive loops on vertices. When visualising connectivity preorders as directed graphs below we always use this format.
\begin{example}\label{ex:q_connectivity_problem_solved}
	We revisit Example \ref{ex:q_connectivity_problem}. By choosing \((0,\widehat{d_0},\widehat{d_1})\)-connection, the associated preorders are, respectively, as shown below.
	\begin{center}
		\begin{tikzpicture}
		\tikzstyle{point}=[circle,thick,draw=black,fill=black,inner sep=0pt,minimum width=2pt,minimum height=2pt]
		\tikzstyle{arc}=[shorten >= 12pt,shorten <= 12pt,->, thick]
		
		\node[] at (0,1.5) {(0)}; 
		\node[] at (-1.2,-0.3) {(1)};
		\node[] at (1.2,-0.3) {(2)};
		\node[] at (-1.2,1) {(01)};
		\node[] at (1.2,1) {(20)};
		\node[] at (0,-1) {(12)};

		\draw[arc] (0,1.5) to (1.2,1);
		\draw[arc] (0,1.5) to (-1.2,1);
		\draw[arc] (-1.2,-0.3) to (-1.2,1);
		\draw[arc] (-1.2,-0.3) to (0,-1);
		\draw[arc] (0,-1) to (1.2,1);
		\draw[arc] (1.2,1) to (-1.2,1);
		\draw[arc] (-1.2,1) to (0,-1);
		\draw[arc] (1.2,-0.3) to (1.2,1);
		\draw[arc] (1.2,-0.3) to (0,-1);
		
		\node[] at (5,1.5) {(a)};
		\node[] at (4,0.25) {(da)};
		\node[] at (5,0.25) {(ba)};
		\node[] at (6,0.25) {(ca)};		
		\node[] at (4,-1) {(d)};
		\node[] at (5,-1) {(b)};
		\node[] at (6,-1) {(c)};
		
		\draw[arc] (5,1.5) to (4,0.25);	
		\draw[arc] (5,1.5) to (5,0.25);	
		\draw[arc] (5,1.5) to (6,0.25);
		\draw[arc] (4,-1) to (4,0.25);
		\draw[arc] (5,-1) to (5,0.25);
		\draw[arc] (6,-1) to (6,0.25);	
		\end{tikzpicture}
	\end{center}
	The \(q\)-connectivity preorders distinguish the complexes. The directionality of the 1-simplices alters the connectivity structure, which is not visible through standard \(Q\)-analysis. In the case of \((0,\widehat{d_0},\widehat{d_0})\)-connectivity we obtain the below preorders.
	\begin{center}
		\begin{tikzpicture}
		\tikzstyle{point}=[circle,thick,draw=black,fill=black,inner sep=0pt,minimum width=2pt,minimum height=2pt]
		\tikzstyle{arc}=[shorten >= 10pt,shorten <= 10pt,->, thick]
		
		\coordinate (0) at (0,1.5);
		\coordinate (1) at (-1.2,-0.3);
		\coordinate (2) at (1.2,-0.3);
		\coordinate (01) at (-1.2,1);
		\coordinate (20) at (1.2,1);
		\coordinate (12) at (0,-1);
		
		\node[] at (0) {(0)}; 
		\node[] at (1) {(1)};
		\node[] at (2) {(2)};
		\node[] at (01) {(01)};
		\node[] at (20) {(20)};
		\node[] at (12) {(12)};

		\draw[arc] (0) to (20);
		\draw[arc] (0) to (01);
		\draw[arc] (1) to (01);
		\draw[arc] (1) to (12);
		\draw[arc] (2) to (20);
		\draw[arc] (1) to (12);
		\draw[arc] (2) to (12);
		
		\coordinate (a) at (5,0.4);
		\coordinate (b) at (3,1.65);
		\coordinate (c) at (7,1.65);
		\coordinate (d) at (5,-1.85);
		\coordinate (da) at (5,-0.75);
		\coordinate (ba) at (4,1.25);
		\coordinate (ca) at (6,1.25);
		
		\node[] at (a) {(a)};
		\node[] at (da) {(da)};
		\node[] at (ba) {(ba)};
		\node[] at (ca) {(ca)};		
		\node[] at (d) {(d)};
		\node[] at (b) {(b)};
		\node[] at (c) {(c)};
		
		\draw[arc] (a) to (da);	
		\draw[arc] (a) to (ba);	
		\draw[arc] (a) to (ca);
		\draw[arc] (d) to (da);
		\draw[arc] (b) to (ba);
		\draw[arc] (c) to (ca);
		\draw[arc,bend left=30] (da) to (ba);
		\draw[arc] (ba) to (da);
		\draw[arc,bend left=30] (ba) to (ca);
		\draw[arc] (ca) to (ba);
		\draw[arc] (da) to (ca);
		\draw[arc,bend left=30] (ca) to (da);
		\end{tikzpicture}
	\end{center}
\end{example}

\begin{example}\label{ex:2_spheres}
	The directed graphs below span 2-dimensional directed flag complexes, each of which has the homotopy type of a 2-sphere. Homology therefore cannot distinguish these complexes. Standard \(Q\)-analysis also sees the undirected simplicial connectivity structures as identical. 	
	\begin{center}
		\begin{tikzpicture}
		\tikzstyle{point}=[circle,thick,draw=black,fill=black,inner sep=0pt,minimum width=2pt,minimum height=2pt]
		\tikzstyle{arc}=[shorten >= 10pt,shorten <= 10pt,->, thick]
		
		\coordinate (0) at (0,1.25);
		\coordinate (3) at (0,-1.25);
		\coordinate (1) at (-1,0);
		\coordinate (2) at (1,0);		
		
		\node[] at (0) {0};
		\node[] at (3) {3};
		\node[] at (1) {1};
		\node[] at (2) {2};
		
		\draw[arc] (0) to (1);
		\draw[arc] (0) to (2);
		\draw[arc] (1) to (3);
		\draw[arc] (2) to (3);
		\draw[arc,bend right=30] (2) to (1);
		\draw[arc,bend right=30] (1) to (2);
		
		\coordinate (a') at (5,1.25);
		\coordinate (b') at (5,-1.25);
		\coordinate (c') at (4,0);
		\coordinate (d') at (6,0);
		
		\node[] at (a') {N};
		\node[] at (b') {S};
		\node[] at (c') {W};
		\node[] at (d') {E};
		
		\draw[arc] (a') to (c');
		\draw[arc] (a') to (d');
		\draw[arc] (b') to (c');
		\draw[arc] (b') to (d');
		\draw[arc,bend right=30] (d') to (c');
		\draw[arc,bend right=30] (c') to (d');
		
		\end{tikzpicture}
	\end{center}
	The digraphs and the complexes are, however, evidently distinct. On the left it is possible to say that there would be two 2-simplicial flow paths from 0 to 3 through simplices (012) and (123), and simplices (021) and (213). On the right these flows are obstructed, but there are circular flows on the upper and lower hemispheres.
	
	The preorders of \((1,\widehat{d_0},\widehat{d_2})\)-connections below show the 2-simplicial flows from 0 to 3 in the left digraph:
	\begin{center}
		\begin{tikzpicture}
		\tikzstyle{point}=[circle,thick,draw=black,fill=black,inner sep=0pt,minimum width=2pt,minimum height=2pt]
		\tikzstyle{arc}=[shorten >= 12pt,shorten <= 12pt,->, thick]
		
		\coordinate (01) at (1,1.75);
		\coordinate (02) at (1,1);
		\coordinate (13) at (1,-2.5);
		\coordinate (23) at (1,-3.25);
		\coordinate (12) at (-1.25,-0.75);
		\coordinate (21) at (3.25,-0.75);
		\coordinate (012) at (0,0);
		\coordinate (021) at (2,0);
		\coordinate (123) at (0,-1.5);
		\coordinate (213) at (2,-1.5);		
		
		\node[] at (01) {(01)};
		\node[] at (02) {(02)};
		\node[] at (13) {(13)};
		\node[] at (23) {(23)};
		\node[] at (12) {(12)};
		\node[] at (21) {(21)};
		\node[] at (012) {(012)};
		\node[] at (021) {(021)};
		\node[] at (123) {(123)};
		\node[] at (213) {(213)};
		
		\draw[arc] (012) to (123);
		\draw[arc] (021) to (213);
		\draw[arc] (02) to (012);
		\draw[arc] (02) to (021);
		\draw[arc] (12) to (012);
		\draw[arc] (12) to (123);	
		\draw[arc] (21) to (021);
		\draw[arc] (21) to (213);
		\draw[arc] (13) to (123);
		\draw[arc] (13) to (213);
		\draw[arc,bend right=20] (01) to (012);
		\draw[arc,bend left=20] (01) to (021);
		\draw[arc,bend left=20] (23) to (123);
		\draw[arc,bend right=20] (23) to (213);
		
		\coordinate (nw) at (7.5,1.75);
		\coordinate (ne) at (7.5,1);
		\coordinate (sw) at (7.5,-2.5);
		\coordinate (se) at (7.5,-3.25);
		\coordinate (we) at (5.25,-0.75);
		\coordinate (ew) at (9.75,-0.75);
		\coordinate (nwe) at (6.5,0);
		\coordinate (new) at (8.5,0);
		\coordinate (swe) at (6.5,-1.5);
		\coordinate (sew) at (8.5,-1.5);
		
		\node[] at (nw) {(NW)};
		\node[] at (ne) {(NE)};
		\node[] at (sw) {(SW)};
		\node[] at (se) {(SE)};
		\node[] at (we) {(WE)};
		\node[] at (ew) {(EW)};
		\node[] at (nwe) {(NWE)};
		\node[] at (new) {(NEW)};
		\node[] at (swe) {(SWE)};
		\node[] at (sew) {(SEW)};
		
		\draw[arc] (ne) to (nwe);
		\draw[arc] (ne) to (new);
		\draw[arc] (we) to (nwe);
		\draw[arc] (we) to (swe);
		\draw[arc] (ew) to (new);
		\draw[arc] (ew) to (sew);	
		\draw[arc] (sw) to (swe);
		\draw[arc] (sw) to (sew);
		\draw[arc,bend right=20] (nw) to (nwe);
		\draw[arc,bend left=20] (nw) to (new);
		\draw[arc,bend left=20] (se) to (swe);
		\draw[arc,bend right=20] (se) to (sew);
		\end{tikzpicture}
	\end{center}
	On the other hand, the \((1,\widehat{d_1},\widehat{d_2})\)-connections reveal the circular flows on the hemispheres:
	\begin{center}
		\begin{tikzpicture}
		\tikzstyle{point}=[circle,thick,draw=black,fill=black,inner sep=0pt,minimum width=2pt,minimum height=2pt]
		\tikzstyle{arc}=[shorten >= 12pt,shorten <= 12pt,->, thick]
		
		\coordinate (01) at (1,1.75);
		\coordinate (02) at (1,1);
		\coordinate (13) at (1,-2.5);
		\coordinate (23) at (1,-3.25);
		\coordinate (12) at (-1.25,-0.75);
		\coordinate (21) at (3.25,-0.75);
		\coordinate (012) at (0,0);
		\coordinate (021) at (2,0);
		\coordinate (123) at (0,-1.5);
		\coordinate (213) at (2,-1.5);		
		
		\node[] at (01) {(01)};
		\node[] at (02) {(02)};
		\node[] at (13) {(13)};
		\node[] at (23) {(23)};
		\node[] at (12) {(12)};
		\node[] at (21) {(21)};
		\node[] at (012) {(012)};
		\node[] at (021) {(021)};
		\node[] at (123) {(123)};
		\node[] at (213) {(213)};
		
		\draw[arc] (012) to (123);
		\draw[arc] (021) to (213);
		\draw[arc] (02) to (012);
		\draw[arc] (02) to (021);
		\draw[arc] (12) to (012);
		\draw[arc] (12) to (123);	
		\draw[arc] (21) to (021);
		\draw[arc] (21) to (213);
		\draw[arc] (13) to (123);
		\draw[arc] (13) to (213);
		\draw[arc,bend right=20] (01) to (012);
		\draw[arc,bend left=20] (01) to (021);
		\draw[arc,bend left=20] (23) to (123);
		\draw[arc,bend right=20] (23) to (213);
		\draw[arc,bend right=20] (021) to (012);
		\draw[arc,bend right=20] (012) to (021);
		
		\coordinate (nw) at (7.5,1.75);
		\coordinate (ne) at (7.5,1);
		\coordinate (sw) at (7.5,-2.5);
		\coordinate (se) at (7.5,-3.25);
		\coordinate (we) at (5.25,-0.75);
		\coordinate (ew) at (9.75,-0.75);
		\coordinate (nwe) at (6.5,0);
		\coordinate (new) at (8.5,0);
		\coordinate (swe) at (6.5,-1.5);
		\coordinate (sew) at (8.5,-1.5);
		
		\node[] at (nw) {(NW)};
		\node[] at (ne) {(NE)};
		\node[] at (sw) {(SW)};
		\node[] at (se) {(SE)};
		\node[] at (we) {(WE)};
		\node[] at (ew) {(EW)};
		\node[] at (nwe) {(NWE)};
		\node[] at (new) {(NEW)};
		\node[] at (swe) {(SWE)};
		\node[] at (sew) {(SEW)};
		
		\draw[arc] (ne) to (nwe);
		\draw[arc] (ne) to (new);
		\draw[arc] (we) to (nwe);
		\draw[arc] (we) to (swe);
		\draw[arc] (ew) to (new);
		\draw[arc] (ew) to (sew);	
		\draw[arc] (sw) to (swe);
		\draw[arc] (sw) to (sew);
		\draw[arc,bend right=20] (nw) to (nwe);
		\draw[arc,bend left=20] (nw) to (new);
		\draw[arc,bend left=20] (se) to (swe);
		\draw[arc,bend right=20] (se) to (sew);
		
		\draw[shorten >= 18pt,shorten <= 18pt,->, thick,bend right=20] (new) to (nwe);
		\draw[shorten >= 18pt,shorten <= 18pt,->, thick,bend right=20] (nwe) to (new);
		\draw[shorten >= 18pt,shorten <= 18pt,->, thick,bend right=20] (sew) to (swe);
		\draw[shorten >= 18pt,shorten <= 18pt,->, thick,bend right=20] (swe) to (sew);
		\end{tikzpicture}
	\end{center}
	The \((1,\widehat{d_0},\widehat{d_1})\)-connection would show the circular flow on the lower hemisphere of the left digraph, whereas the right digraph would be the same as for the \((1,\widehat{d_0},\widehat{d_2})\)-connection above.
\end{example}

The directed \(q\)-connection now adds a new relation on the sets of simplices of ordered simplicial complexes:
\begin{center}
	\renewcommand{\arraystretch}{1.5}
	\begin{tabular}{ c | p{3.5cm} }
		\textbf{Connection} & \textbf{Relation type}\\ 
		\hline
		\((q,\widehat{d_i},\widehat{d_j})\)-connection & preorder \\ \hline
		face poset & partial order \\ 
		\hline
		\(q\)-connection & equivalence relation		
	\end{tabular}
\end{center}
By Definition \ref{def:q_connection_directed} inclusions are always near for any pair \((\widehat{d_i},\widehat{d_j})\). We therefore have a result analogous to Proposition \ref{prop:face_poset_inside_q_connection} saying that the homotopy type of a complex can be reconstructed from its directed \(q\)-connection (see \cite{Bjorner_CW_posets} for the fact that an ordered simplicial complex and its face poset are homotopy equivalent). As the directions \((\widehat{d_i},\widehat{d_j})\) play no role in \(q\)-nearness by inclusion, the proof is exactly as in Proposition \ref{prop:face_poset_inside_q_connection}.
\begin{proposition}
	Let \(K\) be an ordered simplicial complex and \(F\) its face poset. Let \(R_q\) denote the relation of being \(q\)-near along any \((\widehat{d_i},\widehat{d_j})\) on the set of simplices \(K_q\). Then \(F \subset \bigcup_{q=0}^{\dim(K)} R_q\).
\end{proposition}
The relations in the table above can now be seen to give a hierarchy of relations for collections of simplices. The face poset can be reconstructed from the \((q,\widehat{d_i},\widehat{d_j})\)-connections, which again can be obtained from the equivalence relation of standard \(q\)-connection by restrictions given by \((\widehat{d_i},\widehat{d_j})\). 
\begin{remark}
	Any simplicial complex can be turned into an ordered simplicial complex by fixing a linear ordering on the vertex set. Different orderings yield isomorphic ordered complexes. Thus from the point of view of classical \(Q\)-analysis these are all indistinguishable. The directed \(Q\)-analysis on the other hand can make a difference between different orderings. We have not explored whether the directed \(Q\)-analysis applied to different orderings of simplicial complexes might be interesting; for us the interest rises from the explicit orderings given by underlying digraphs.
\end{remark}

The following is helpful in understanding aspects of the \(q\)-connections for different choices of \((\widehat{d_i},\widehat{d_j})\).

\begin{lemma}\label{lem:q_preorders_structure}
	Let \(\sigma,\tau \in K_q\) and assume \(\sigma \not\hookrightarrow \tau\) and \(\tau \not\hookrightarrow \sigma\). Then we have the following properties.
	\begin{enumerate}
		\item If the ordered pair \((\sigma,\tau)\) is \((q,\widehat{d_i},\widehat{d_j})\)-near, then the ordered pair \((\tau,\sigma)\) is \((q,\widehat{d_j},\widehat{d_i})\)-near.
		\item If the ordered pair \((\sigma,\tau)\) is \((q,\widehat{d_i},\widehat{d_i})\)-near, then so is \((\tau,\sigma)\).
	\end{enumerate}
\end{lemma}
\begin{proof}
	For 1., if \(\widehat{d_i}(\sigma) \hookleftarrow \alpha \hookrightarrow \widehat{d_j}(\tau),\) for some \(q\)-simplex \(\alpha \in K\), then it is immediate that \(\widehat{d_j}(\tau) \hookleftarrow \alpha \hookrightarrow \widehat{d_i}(\sigma)\). Property 2. follows likewise from the symmetry of the relation \(\widehat{d_i}(\sigma) \hookleftarrow \alpha \hookrightarrow \widehat{d_i}(\tau)\).
\end{proof}

We finish this section by defining directed analogs of pseudomanifolds. As mentioned in the introduction, one motivation is the hypothesis that neural dynamics on a directed network resides on a low-dimensional manifold, and its manifestation on the network level might be of interest.
\begin{definition}
	An ordered simplicial complex \(K\) is an \(n\)\textbf{-pseudomanifold along \((\widehat{d_i},\widehat{d_j})\)} if all the maximal simplices are \(n\)-simplices, each \((n-1)\)-simplex is a face of exactly two \(n\)-simplices, and any two \(n\)-simplices are \((n-1,\widehat{d_i},\widehat{d_j})\)-connected.
\end{definition}

\begin{definition}
	An ordered simplicial complex \(K\) is an \(n\)\textbf{-pseudomanifold along \((\widehat{d_i},\widehat{d_j})\) with boundary} if all the maximal simplices are \(n\)-simplices, each \((n-1)\)-simplex is a face of at most two \(n\)-simplices, and any two \(n\)-simplices are \((n-1,\widehat{d_i},\widehat{d_j})\)-connected.
\end{definition}

\begin{example}
	The directed flag complex below is a 2-pseudomanifold along \((\widehat{d_0},\widehat{d_1})\) with boundary.
	\begin{figure}[h!]\label{fig:example_directed_pseudomanifold}
		\begin{center}
			\begin{tikzpicture}
			\tikzstyle{point}=[circle,thick,draw=black,fill=black,inner sep=0pt,minimum width=2pt,minimum height=2pt]
			\tikzset{dirarrow/.style={postaction={decorate,decoration={markings,mark=at position .7 with {\arrow[#1]{stealth}}}}}}
			\tikzstyle{arc}=[shorten >= 8pt,shorten <= 8pt,->]
			
			\coordinate (a) at (-1.4,-0.1);
			\coordinate (b) at (0,1);
			\coordinate (c) at (1.4,-0.1);
			\coordinate (d) at (0.9,-1.6);
			\coordinate (e) at (-0.9,-1.6);
			\coordinate (f) at (0,-0.4);
			
			\draw[fill=magenta,opacity=0.3] (a) -- (b) -- (c) -- (d) -- (e) --cycle;

			\draw[postaction={dirarrow}] (a) to (b);
			\draw[postaction={dirarrow}] (b) -- (c);
			\draw[postaction={dirarrow}] (c) -- (d);
			\draw[postaction={dirarrow}] (d) -- (e);
			\draw[postaction={dirarrow}] (e) -- (a);
			
			\draw[postaction={dirarrow}] (a) -- (f);
			\draw[postaction={dirarrow}] (b) -- (f);
			\draw[postaction={dirarrow}] (c) -- (f);
			\draw[postaction={dirarrow}] (d) -- (f);
			\draw[postaction={dirarrow}] (e) -- (f);
			
			\node[point] at (a) {};
			\node[point] at (d) {};
			\node[point] at (c) {};
			\node[point] at (b) {};
			\node[point] at (e) {};
			\node[point] at (f) {};
			\end{tikzpicture}
		\end{center}
	\end{figure}
\end{example}
The directed cycle in Example \ref{ex:q_connectivity_problem} is a 1-pseudomanifold along \((\widehat{d_0},\widehat{d_1})\), as shown by its connectivity preorder in Example \ref{ex:q_connectivity_problem_solved}; in fact any directed cycle gives an example of such a pseudomanifold. In general, it seems very difficult to come up with examples of higher-dimensional directed pseudomanifolds. For example, neither of the 2-spheres in Example \ref{ex:2_spheres} is a pseudomanifold along any \((\widehat{d_i},\widehat{d_j})\). It is an interesting open question whether directed \(n\)-pseudomanifolds can be constructed and whether they exist in some naturally occurring digraphs. 

We can say something about the structure of pseudomanifolds. Note that in the case of simplicial complex pseudomanifolds of Section \ref{sec:q_connectivity}, and considering classical \(q\)-connectivity, the flag complexes of the associated \(q\)-graphs can be higher than 1-dimensional. The illustration below shows a simple example of a 2-pseudomanifold with boundary, whose flag complex of the \(1\)-graph contains 2-simplices. In comparison, the \((q,\widehat{d_i},\widehat{d_j})\)-graphs of directed pseudomanifolds have a more restricted structure.
\begin{figure}[h!]\label{fig:pseudomanifold_has_2simplex}
	\begin{center}
		\begin{tikzpicture}
		\tikzstyle{point}=[circle,thick,draw=black,fill=black,inner sep=0pt,minimum width=2pt,minimum height=2pt]
		
		\coordinate (a) at (-1,0);
		\coordinate (b) at (1,0);
		\coordinate (c) at (0,1.45);
		\coordinate (d) at (0,0.6);
		
		\draw[fill=magenta,opacity=0.3] (a) -- (b) -- (c) --cycle; 		
		
		\draw[] (a) -- (b);
		\draw[] (b) -- (c);
		\draw[] (c) -- (a);
		\draw[] (d) -- (c);
		\draw[] (d) -- (a);
		\draw[] (d) -- (b);
		
		\node[point] at (a) {};
		\node[point] at (d) {};
		\node[point] at (c) {};
		\node[point] at (b) {};
		\end{tikzpicture}
	\end{center}
\end{figure}

\begin{proposition}
	Let \(K\) be an \(n\)-pseudomanifold along \((\widehat{d_i},\widehat{d_j})\), with or without boundary. Let \(\GG\) be the subgraph of the \((n-1,\widehat{d_i},\widehat{d_j})\)-graph induced by the \(n\)-simplices. Then the directed flag complex of \(\GG\) is 1-dimensional.
\end{proposition}
\begin{proof}
	We show that \(\GG\) cannot have 2-simplices in its associated directed flag complex. So assume that there is a 2-simplex \((\sigma,\tau,\kappa)\), spanned by the given \(n\)-simplices as vertices. Then there are directed edges \((\sigma,\tau)\), \((\sigma,\kappa)\), and \((\tau,\kappa)\). Note that the face \(\widehat{d_i}(\sigma)\) is a unique \((n-1)\)-simplex \(\alpha\). The \((n-1,\widehat{d_i},\widehat{d_j})\)-nearness inducing the above directed edges then gives us the following relations:
	\[\widehat{d_j}(\kappa) \hookleftarrow \alpha \hookrightarrow \widehat{d_i}(\sigma) \hookleftarrow \alpha \hookrightarrow \widehat{d_j}(\tau).\]
	But this gives a contradiction to the pseudomanifold conditions since the \((n-1)\)-simplex \(\alpha\) is now a face of three \(n\)-simplices.
\end{proof}

\subsection{Directed structure systems}\label{sec:dir_struct_sys}
The basic aim of standard \(Q\)-analysis is to attach a simplicial complex \(K\) with its first structure vector \(\mathbf{Q}(K)\) encoding the connectivity information by the number of \(q\)-connected components. We introduce in this section an analog for directed \((q,\widehat{d_i},\widehat{d_j})\)-connection. As this outputs a collection of preorders, the information is no longer summarised simply by numbers of equivalence classes. In turn, a richer structure arises that sees the simplicial connectivity through a collection of topological spaces arising from the preorders. 

Let \(D=\text{dim}(K)\) for an ordered simplicial complex \(K\). Each choice of \(0 \leq q \leq D-1\) corresponds to \(|\{\widehat{d_i}\}_{i=0}^{D}| \times|\{\widehat{d_j}\}_{j=0}^{D}|\) choices of ordered pairs of directions. When \(q=D\) all the simplices are maximal and therefore only connected to themselves for any pair \((\widehat{d_i},\widehat{d_j})\), and there is essentially only one preorder. Combined there are \(D \times ((D+1)\times(D+1)) +1 = D^3 + 2D^2 + D +1\) possible \((q,\widehat{d_i},\widehat{d_j})\)-connections. We denote by \(Q_K\) the set of all the possible triples \((q,\widehat{d_i},\widehat{d_j})\). We let \textbf{Preorders} stand for the collection of finite preordered sets.

\begin{definition}
	Let \(K\) be an ordered simplicial complex. The \textbf{first structure map} of \(K\) is the association 
	\[\Gamma_K \colon Q_K \ra \textbf{Preorders}.\]
\end{definition}
There is a famous bijection, as first noted by Alexandroff, of preorders and topological spaces with Alexandroff topologies. We recall here basic facts from the theory of finite topological spaces \cite{Barmak_book}, which allows us to see the \((q,\widehat{d_i},\widehat{d_j})\)-connections through topology. In the finite setting we have the following correspondence.
\begin{proposition}\label{prop:topo_preord_bijection}
	Finite preorders are in bijection with finite topological spaces. Moreover, finite partial orders, or posets, are in bijection with finite \(T_0\) topological spaces.
\end{proposition}
When studying finite spaces from a homotopical point of view, the following proposition says that we can restrict to \(T_0\)-spaces \cite{Barmak_book}. The exhibited equivalence relation is a canonical way of turning a preorder into a partial order \cite[Proposition 8.13]{Schroeder_book}. The relation \(x \le y\) in topological terms means that \(x\) is in every open set containing \(y\).
\begin{proposition}\label{prop:Kolmogorov_quotient}
	Let \(P\) be a finite space. Let \(P/\!\!\sim\) be the quotient where \(x \sim y\) if and only if \(x \le y\) and \(y \le x\). Then \(P/\!\!\sim\) is a \(T_0\)-space and the quotient map \(q \colon P \ra P/\!\!\sim\) is a homotopy equivalence.
\end{proposition}
Any preorder can also be seen as a digraph, where we put an edge \((x,y)\) if \(x \le y\). Proposition \ref{prop:topo_preord_bijection} then gives a way to see finite topological spaces as (transitive) digraphs (now with self-loops \((x,x)\)). Finite \(T_0\) spaces are further in correspondence with (transitive) directed acyclic graphs, or DAGs. As noted in \cite[Lemma 6]{Berghammer_Winter_2019}, Proposition \ref{prop:Kolmogorov_quotient} has a counterpart in terms of digraphs, called the condensation of a digraph.
\begin{definition}
	A \textbf{strongly connected component} in a digraph \(\GG\) is an induced subgraph \(\GG'\) such that for every pair of vertices \(x\) and \(y\) in \(\GG'\) there is a path \(x \ra y\) in \(\GG'\).
\end{definition}
The strongly connected components are the equivalence classes of the relation of being strongly connected on the vertices of \(\GG\), i.e.\ having directed paths between any ordered pair of vertices. The ensuing partition then enables to construct the quotient graph without directed cycles.
\begin{definition}\label{def:condensation}
	The \textbf{condensation} \(c(\GG)\) of digraph \(\GG\) has as its vertices the strongly connected components of \(\GG\). Two vertices \(X\) and \(Y\) have a directed edge \((X,Y)\) in \(c(\GG)\) if there is an edge \((x,y)\) in \(\GG\) for some \(x \in X\) and \(y \in Y\). The digraph \(c(\GG)\) is acyclic.
\end{definition}
The crux of the preceding discussion is that the first structure map \(\Gamma_K\) endows the set of ordered simplices with a collection of finite topological spaces. Up to homotopy it is enough to study these spaces in the form of their condensed partial orders. The weak homotopy type of a partial order \(P\) is given by its \textbf{order complex} whose simplices are the chains, or the non-empty totally ordered subsets, of \(P\). By Proposition \ref{prop:face_poset_inside_q_connection} the homotopy type of an ordered simplicial complex \(K\) can be reconstructed from the map \(\Gamma_K\). The \((q,\widehat{d_i},\widehat{d_j})\)-connections, however, exhibit homotopy types different from the original complex.

\begin{example}\label{ex:2_spheres_connection_topology}
	Recall the \((1,\widehat{d_1},\widehat{d_2})\)-connections of the 2-spheres in Example \ref{ex:2_spheres}. Passing to condensations and order complexes of the \((1,\widehat{d_1},\widehat{d_2})\)-preorders, the homotopy type associated to the sphere given by vertices \{0,1,2,3\} is a wedge of circles \(S^1 \vee S^1\), while that of the sphere given by vertices \{N,W,E,S\} is \(S^1\).
\end{example}

Our construction of \((q,\widehat{d_i},\widehat{d_j})\)-connections thus reveals new homotopy types arising from the directionality of the underlying directed graph. By Example \ref{ex:2_spheres_connection_topology} these homotopy types can differ from those seen by the simplicial homology of directed flag complexes, giving us a set of new topological invariants. In analogy to the remark after Theorem \ref{prop:q_connectivity_equivalence} the \((q,\widehat{d_i},\widehat{d_j})\)-connections as preorders give us immediately categories \(\Gamma_{(q,\widehat{d_i},\widehat{d_j})}\).

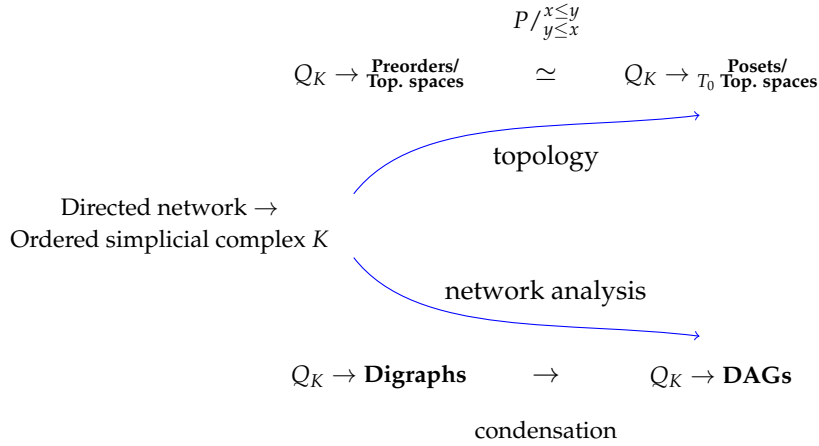
\begin{figure}[h]
	\begin{tikzpicture}
	\tikzstyle{arc}=[shorten >= 8pt,shorten <= 8pt,->]
	\node[align=center] at (-0.25,0) {\small{Directed network \(\ra\)} \\ \small{Ordered simplicial complex \(K\)}};
	\node at (2.5,2) {\small{\(Q_K \ra \substack{\textbf{Preorders/}\\\textbf{Top. spaces}}\)}};
	\node at (2.5,-2) {\small{\(Q_K \ra \textbf{Digraphs}\)}};
	
	\node at (4.7,2) {\(\simeq\)};
	\node at (4.7,2.7) {\small{\(P/\substack{x \le y \\ y \le x}\)}};
	
	\node at (7,2) {\small{\(Q_K \ra \substack{\textbf{Posets/}\\T_0 \textbf{ Top. spaces}}\)}};
	\node at (7,-2) {\small{\(Q_K \ra \textbf{DAGs}\)}};
	
	\node at (4.7,-2) {\(\ra\)};
	\node at (4.7,-2.7) {\small{condensation}};
	
	\draw[arc,blue] (2,0.2) .. controls (3,1.5) and (5,1.25) .. (7,1.5);
	\draw[arc,blue] (2,-0.2) .. controls (3,-1.5) and (5,-1.25) .. (7,-1.5);
	
	\node at (4.7,0.9) {topology};
	\node at (4.7,-0.9) {network analysis};
	\end{tikzpicture}
	\caption{The topological point of view and the digraph point of view offer two complementary ways to study the simplicial \((q,\widehat{d_i},\widehat{d_j})\)-connections of ordered simplicial complexes and ultimately of networks.}
	\label{fig:big_picture}
\end{figure}

Returning to our second motivating question of higher simplicial paths, the viewpoint of the first structure map \(\Gamma_K\) as a collection of digraphs facilitates this with path searches. Here the passage from a \((q,\widehat{d_i},\widehat{d_j})\)-preorder to its acyclic condensation offers considerable algorithmic advantage by removing directed cycles. We exploit this in the network analyses in the next section. The preceding discussion and our view on the \((q,\widehat{d_i},\widehat{d_j})\)-connections are summarised in Figure \ref{fig:big_picture}.

\section{Applications to network analysis}\label{sec:applications}
Our interest in network analysis stems from finding simplicial pathways. As depicted in Figure \ref{fig:example_pathways}, various edge paths of a network are supported on directed simplices, which then induce their own higher dimensional paths. Our construction of \((q,\widehat{d_i},\widehat{d_j})\)-connections, now seen as digraphs of simplices, is appropriate for extracting these paths. Different choices of \(q\), \(i\) and \(j\) allow to emphasise different features of directionality, in the same vein as in the previous section the \((q,\widehat{d_i},\widehat{d_j})\)-preorders or -posets as topological spaces let us see different homotopy types. 

We are particularly inspired by the field of topological neuroscience, and in this section we apply our construction to various brain networks. Specifically we compute the longest simplicial paths in the induced directed flag complexes with some choices of the \((q,\widehat{d_i},\widehat{d_j})\)-connections. The simplicial path analysis is seen to show different connectivity properties between the networks. To gain more insight into the simplicial paths themselves we introduce so called path fraction which reveals interesting information about the tightness of the simplicial connectivity along the paths. For simplified notation, in this section we write \((q,\widehat{d_i},\widehat{d_j})\) as \((q,i,j)\). We remark that the aim of this section is not to make any conclusive neuroscientific claims, but to demonstrate the potential of the developed theory for network science and topological neuroscience.

\subsection{Data description}\label{subsec:data_description}
We used in our analyses the following networks of neurons:

\textbf{Blue Brain reconstruction \cite{BB_cell_paper}.} The Blue Brain microcircuit is a biologically validated digital reconstruction of the neuronal connectivity in a small volume of a rat somatosensory cortex. The network spans across 6 neuronal layers, each with its characteristic neuron types and connectivity patterns. The topological analysis of the Blue Brain reconstruction in terms of simplicial homology has previously been done in \cite{Frontiers_paper}. The microcircuit also facilitates simulation of neuronal activity. Recently, various simulated activities were classified with high accuracy using feature vectors constructed from the network structure \cite{Tribes_math,Tribes_neuroscience}. The microcircuit is obtainable from \cite{BB_portal}. The graph we used has 31,346 vertices and 7,803,528 edges.

\textbf{C. elegans \cite{C_elegans}.} The neuronal network of the worm Caenorhabditis elegans is sufficiently small and has been reconstructed as a whole. Moreover, the neurons and their synaptic connectivity are to a high degree consistent across individuals. A particular feature of the network is an overrepresentation of reciprocally connected triangle motifs and an underrepresentation of directed triangle cycles. The synaptic connectivity is obtainable from \cite{WormAtlas} and we used the steps in \cite{Govc_2020} to construct the directed graph. The graph we used has 279 vertices and 2,194 edges.

\textbf{Allen mouse \cite{Allen_mouse}.} The Allen Institute's model of the mouse primary visual cortex on area V1 contains \(\sim\)230,000 neurons. The model enables simulating cortical behaviour with arbitrary visual stimuli. The network model can be constructed by employing the Brain Modeling ToolKit (BMTK) \cite{BMTK} following the instructions in \cite{Allen_mouse}. We used a slightly modified construction of \cite{simplex_closing} to allow sampling a more feasible size subnetwork. The graph we used has 69,335 vertices and 6,326,900 edges.

\textbf{Drosophila \cite{drosophila}.} This network is a reconstruction of a portion of the central brain of the fruit fly Drosophila melanogaster. The connectivity graph consists of approximately 25,000 vertices and 3 million edges. Similar to the C. elegans graph, the Drosophila has an overrepresentation of reciprocal connections, as well as a significant occurrence of directed cliques compared to a random graph with similar connection probability. The network data is available at \cite{Janelia} and we used the version v1.2 of the reconstruction. The Drosophila network presented major computational challenges. Due to its density and large number of cliques there was a proliferation of simplices amounting to such memory requirements that we were not able to compute them for the full graph. We found feasible simplex counts for the subsequent connectivity analysis by restricting to a subgraph induced by vertices of total degree at most 70. The graph we used has 3,243 vertices and 16,358 edges.

\begin{figure}
	\includegraphics[width=\textwidth]{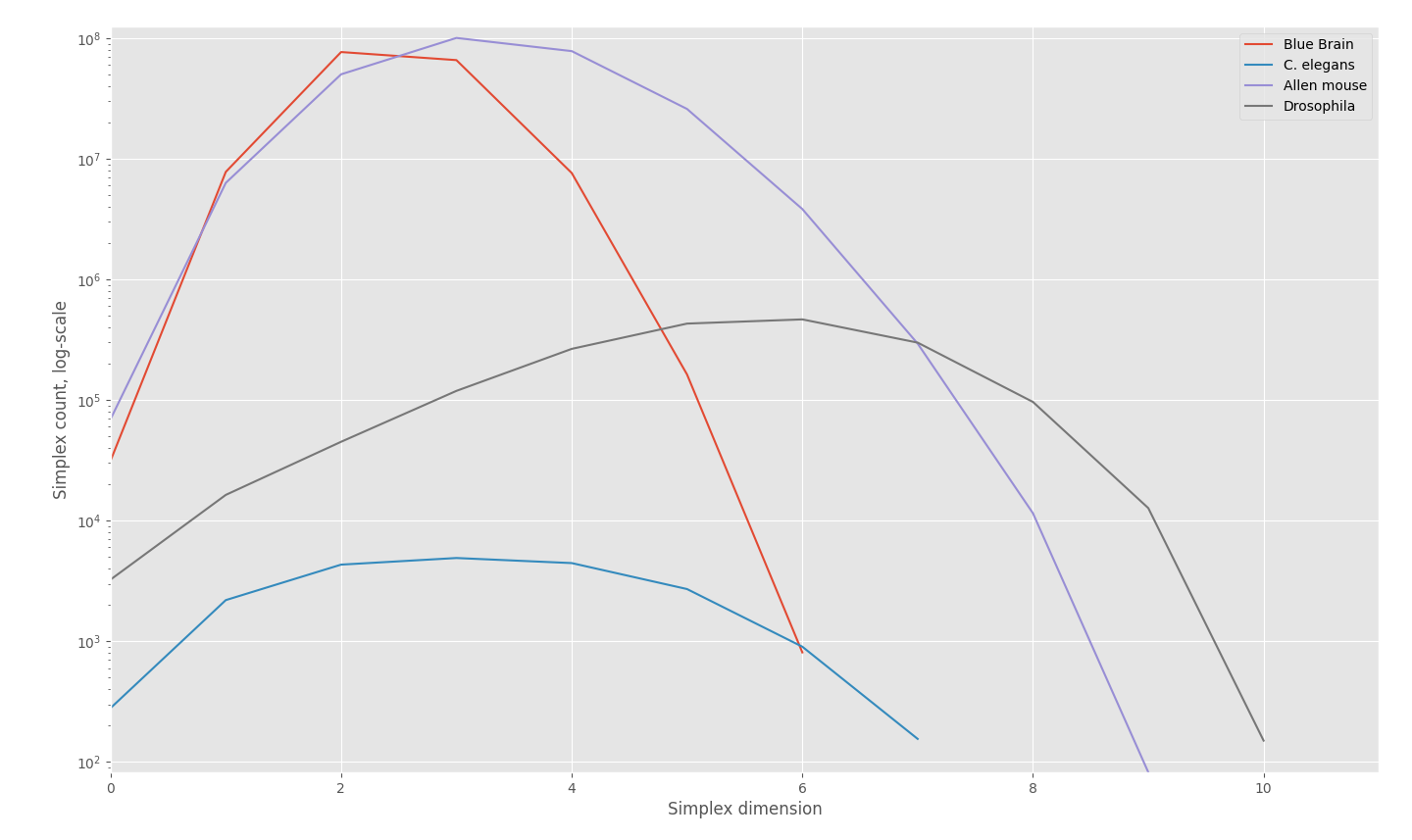}
	\caption{Simplex counts in each dimension for the analysed networks. The \(y\)-axis is on log-scale.}
	\label{fig:simplex_counts}
\end{figure}

Each network was a simple directed graph without self-loops and represented by its adjacency matrix. For creating the directed flag complexes we used the Flagser-count software \cite{flagser_count} that allows outputting and storing the directed simplices needed for the \((q,i,j)\)-connectivity analysis. We plot the simplex counts in each dimension in Figure \ref{fig:simplex_counts}.

\subsection{Computing simplicial paths}
Finding the longest path in a general directed graph is an NP-hard problem. However, for a directed acyclic graph (DAG) the longest path can be found in linear time in terms of the number of vertices and edges in the graph. We therefore employ the condensation \(c(\GG)\) (Definition \ref{def:condensation}) which produces a DAG (partial order) from a \((q,i,j)\)-connectivity digraph \(\GG\) (preorder). We are always working with a Hasse diagram form so that the connectivity digraphs only contain the edges for simplices that are \((q,i,j)\)-near. 

We used Algorithm \ref{alg:longest_path} for computing the simplicial paths. Essentially it traverses the longest path \(P\) in \(c(\GG)\), and forms an augmented path \(\bar{P}\) by finding for every node \(p \in P\) the longest of the shortest paths through \(p\) (every node of \(c(\GG)\) consists of vertices of \(\GG\) and therefore induces a subgraph of \(\GG\)). Note that our longest of the shortest paths within \(p\) differs from the one realising the diameter of \(p\) (graph's diameter is the longest shortest path between \emph{any} two vertices) in that the shortest paths are only computed for those pairs of vertices which enforce continuity of the simplicial path along \(P\), i.e.\ \((q,i,j)\)-nearness for consecutive simplices. By successors(\(v\)) we denote those vertices \(w\) such that \((v,w)\) is an edge in \(\GG\). We implemented Algorithm \ref{alg:longest_path} in Python using the NetworkX library. Note that the graph theoretical length of the longest path \(P\) corresponds to the height, or the inductive dimension \cite{Berghammer_Winter_2019}, of the associated topological space. The augmented path \(\bar{P}\) then realises the height as an actual connected path of simplices.
\begin{algorithm}
	\caption{Simplicial path in a \((q,i,j)\)-digraph}\label{alg:longest_path}
	\begin{algorithmic}
		\State $\GG \gets (q,i,j)$-digraph
		\State $c(\GG) \gets $ condensation of \(\GG\)
		\State $P \gets $ the longest path in \(c(\GG)\)
		\State \(L \gets \) initialized empty list \Comment{List for simplices on the simplicial path}
		\State \(k \gets 0\) \Comment{For tracking position along \(P\)}
		\For{node $p \in P$}
		\If{$|p|=1$}
		
		\State append \(L\) with \(p\)
		\State current\textunderscore node \( \gets p\)
		
		\ElsIf{\(|p| > 1\) and \(k=0\)}
		\State targets \( \gets \) \{nodes \(v \in p \ | \ \text{successors}(v) \cap P[1] \ne \emptyset\)\}
		\State \(\GG_p \gets\) subgraph of \(\GG\) induced by vertices in \(p\)
		\State \(L_p \gets\) \(\underset{length}{\mathrm{max}}\)\{shortest paths in \(\GG_p\) from \(v \in p\) to \(w \in\) targets\}
		\State append \(L\) with \(L_p\)
		\State current\textunderscore node \( \gets L_p[-1]\) \Comment{Last element in the list \(L_p\)}
		
		\ElsIf{\(|p| > 1\) and \(0 < k < length(P)-1\)}
		\State sources \(\gets\) successors(current\textunderscore node) \(\cap\) \(p\) 
		\State targets \( \gets \) \{nodes \(v \in p \ | \ \text{successors}(v) \cap P[k+1] \ne \emptyset\)\}		 
		\State \(\GG_p \gets\) subgraph of \(\GG\) induced by vertices in \(p\)
		\State \(L_p \gets\) \(\underset{length}{\mathrm{max}}\)\{shortest paths in \(\GG_p\) from \(v \in \) sources to \(w \in\) targets\}
		\State append \(L\) with \(L_p\)
		\State current\textunderscore node \( \gets L_p[-1]\)
		
		\ElsIf{\(|p| > 1\) and \(k = length(P)-1\)}
		\State sources \(\gets\) successors(current\textunderscore node) \(\cap\) \(p\) 	 
		\State \(\GG_p \gets\) subgraph of \(\GG\) induced by vertices in \(p\)
		\State \(L_p \gets\) \(\underset{length}{\mathrm{max}}\)\{shortest paths in \(\GG_p\) from \(v \in \) sources to \(w \in p\)\}
		\State append \(L\) with \(L_p\)
		
		\EndIf
		
		\(k \gets k + 1\)
		\EndFor
	\end{algorithmic}
\end{algorithm}

After computing the simplicial path with Algorithm \ref{alg:longest_path}, we wish to gain more insight about the connectivity of simplices along the path \(\bar{P}\). To this end we define the following measure capturing the number of distinct vertices within the simplices on the path relative to the number of vertices in an ideal path of simplices with dimensions the same as in \(\bar{P}\), but where the simplices are connected strictly through \(q\)-faces. Lower path fraction is thus an indication of stronger connectivity between simplices.
\begin{definition}
	Let \((\sigma_0,\sigma_1,\dots,\sigma_n)\) be a simplicial path with respect to \((q,i,j)\)-connection, and let \(V = \bigcup_{i=0}^n \text{Vert}(\sigma_i)\) be the set of distinct vertices on the path. Define 
	\[
	s_i = \begin{cases}
	\dim(\sigma_i)+1, \text{ if } i = 0;\\
	\dim(\sigma_{i+1})-\dim(\sigma_i), \text{ if } \sigma_i \hookrightarrow \sigma_{i+1};\\
	\dim(\sigma_i) - q, \text{ otherwise.}
	\end{cases}\]
	Then the \textbf{path fraction} is given by
	\[\frac{|V|}{\sum_{i=0}^{n} s_i}.\]
\end{definition}

\subsection{Results}
We computed the longest simplicial paths for various \((q,i,j)\)-connections. The results for path lengths are shown in Figure \ref{fig:results_lens} and for path fractions in Figure \ref{fig:results_pfs}; more detailed analysis is given below. It seems more interesting to study simplicial paths with respect to higher values of \(q\) due to the sparsity of higher-dimensional simplices; as shown in Figure \ref{fig:simplex_counts} the simplex counts tend to peak at smaller dimensions. Higher \(q\) also requires higher-dimensional connecting faces; finding long simplicial paths would therefore indicate that the network structure supports certain type of clustering on the simplicial level. We also know from Proposition \ref{prop:directed_q_connectivity_basic_properties} that being \(q\)-connected implies \(p\)-connectivity for any smaller \(p\). Higher \(q\) makes the computation time for the connectivity preorders also feasible due to smaller number of simplices. Our choices for \(q\) are detailed below. Both \(i\) and \(j\) ranged from 0 to 5. We excluded the case \(i=j\) due to the symmetry property in Lemma \ref{lem:q_preorders_structure}, since the induced highly cyclic structure of the connectivity preorders might result in less meaningful simplicial paths.

\begin{figure}
	\includegraphics[width=\textwidth,trim=1.7cm 0cm 1.7cm 0cm,clip]{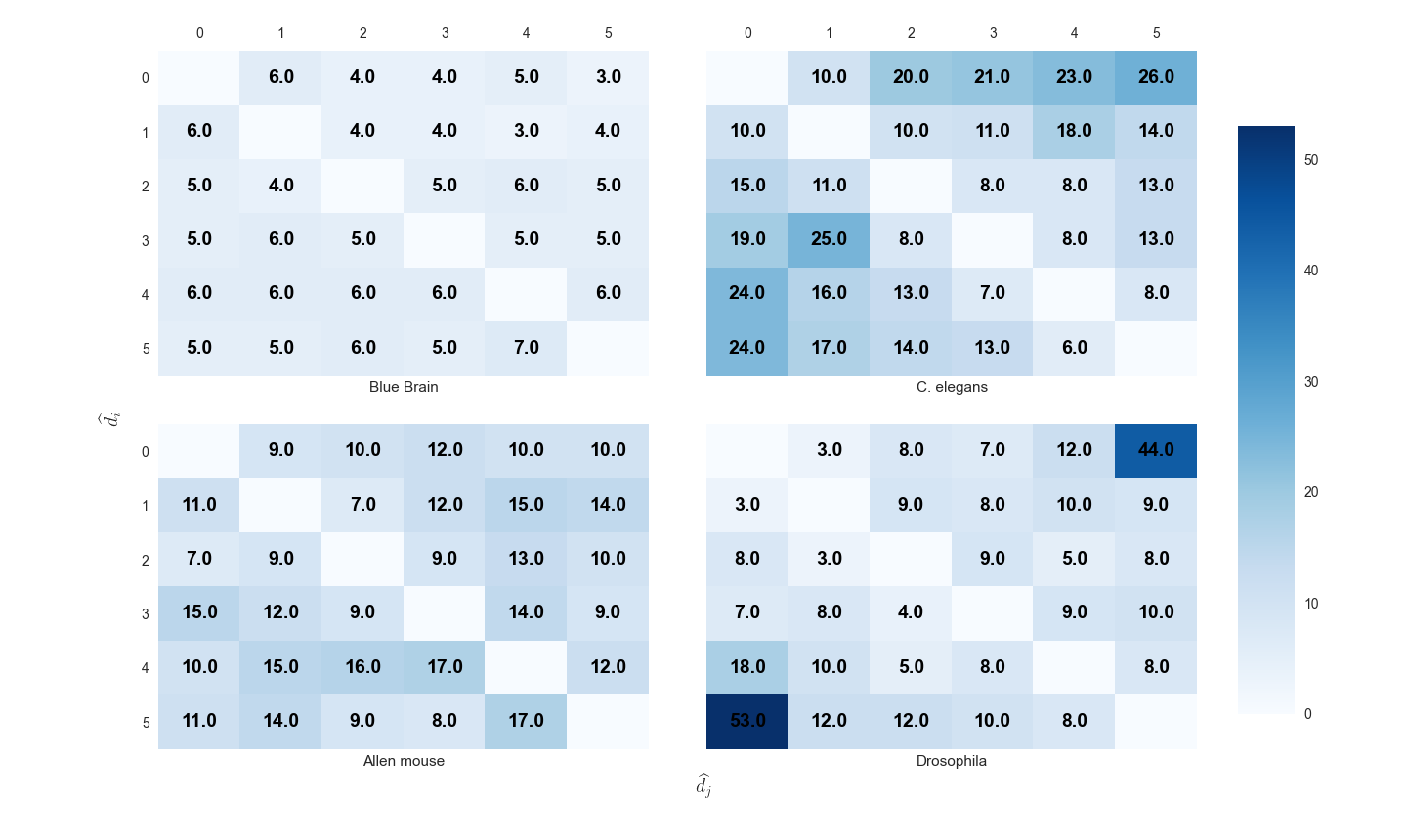}
	\caption{Simplicial path lengths for all analysed networks. The \(\widehat{d_i}\) labels are the same for all rows with values shown on the left; \(\widehat{d_j}\) labels are the same for all columns with values on top. The values of \(q\) used was 4 for Blue Brain and C. elegans, 7 for Allen mouse, and 8 for Drosophila.}
	\label{fig:results_lens}
\end{figure} 
\begin{figure}
	\includegraphics[width=\textwidth,trim=1.7cm 0cm 1.7cm 0cm,clip]{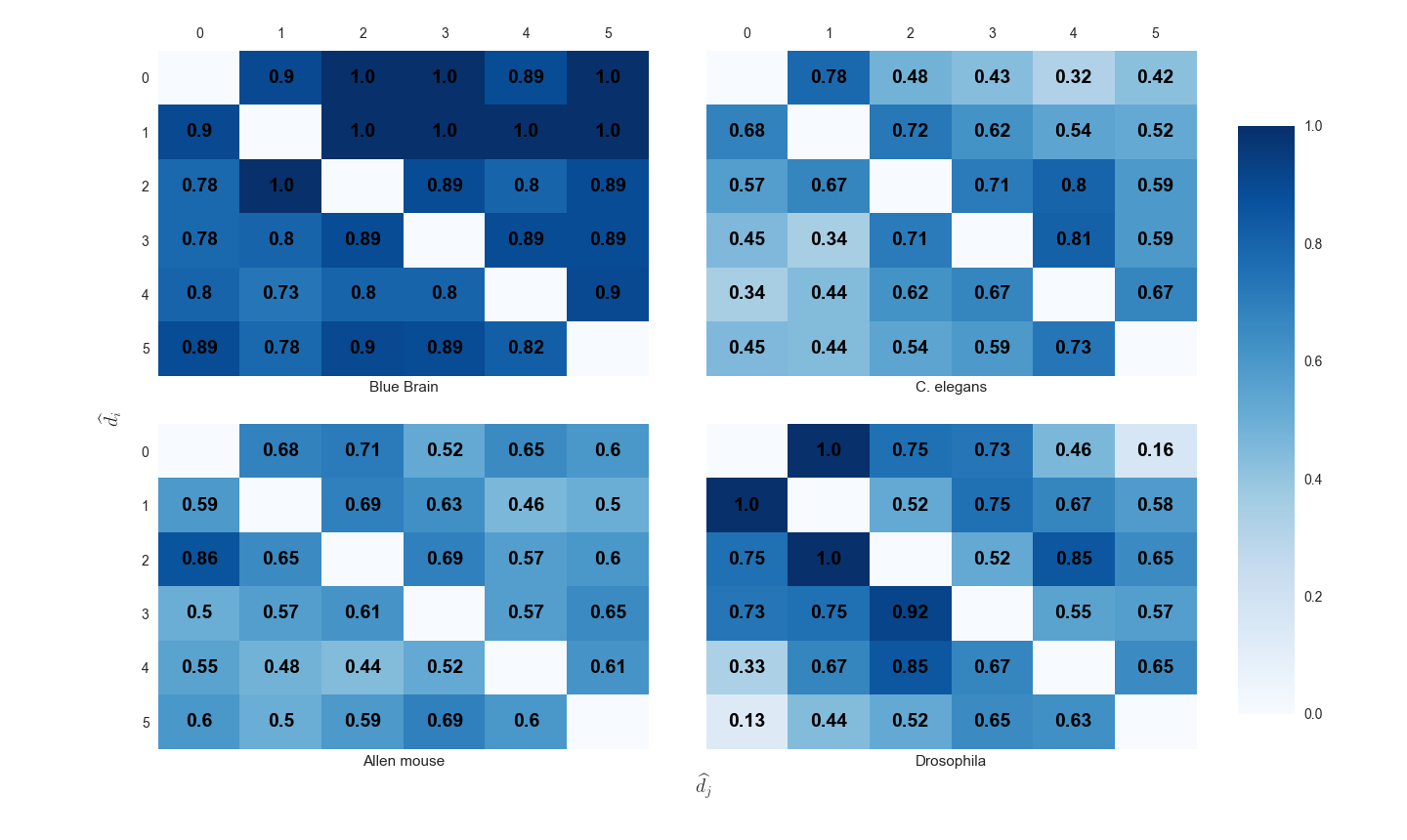}
	\caption{Path fractions for all analysed networks. The \(\widehat{d_i}\) and \(\widehat{d_j}\) labels are shared between rows and columns as in Figure \ref{fig:results_lens}.}
	\label{fig:results_pfs}
\end{figure} 

\textbf{Blue Brain.} As noted in Section \ref{subsec:data_description}, the layers in the Blue Brain network are characterised by their neuron types and their connectivity. The number of simplices per layer also increases from Layer 1 to Layer 6. For a reasonable number of simplices we restricted to Layers 1-4, with \(q=4\). 

The Blue Brain network is characterised by very short path lengths. As shown in Figure \ref{fig:results_lens}, the path lengths are also rather homogeneous over different directions \(\widehat{d_i}\) and \(\widehat{d_j}\). Path fraction analysis in Figure \ref{fig:results_pfs} shows high values, with many exactly 1; the paths are hence close to ideal paths in terms of \(q\)-connectivity. The results indicate that the Blue Brain network is rather void of strong simplicial connectivity and clustering.

We also analysed the (4,0,5)- and (5,0,6)-connectivities with respect to deepening layers: first for the subnetwork induced by Layer 2, then for Layers 2 and 3, then for Layers 2-4 etc. Note that we skipped Layer 1 since it contains no simplices beyond dimension 2. The aim was to investigate whether the layers show any change in the connectivity structure. The choices of \(i\) and \(j\) are motivated by the simplicial source-to-sink paths, an example of which is the \((1,0,2)\)-connection on the left side of Figure \ref{fig:example_pathways}: for any ordered pair of 2-simplices \((\sigma,\tau)\) we take from \(\sigma\) the face in the direction of the source and from \(\tau\) the face in the direction of the sink. The results are shown in Table \ref{table:BlueBrain}. The path lengths are again negligible. Of course more comprehensive analysis with different directions would be needed for the layers as well.
\begin{table}[h!]
	\centering
	\begin{tabular}{|c|c|c|c|c|c|} 
		\hline
		& Layer 2 & Layer 2-3 & Layer 2-4 & Layer 2-5 & Layer 2-6 \\
		\hline
		(4,0,5) & 3 & 3 & 3 & - & - \\ 
		\hline
		(5,0,6) & 1 & 1 & 2 & 2 & 3 \\ 
		\hline
	\end{tabular}
	\caption{Simplicial path lengths in the Blue Brain network with respect to deepening layers. The symbol - means that we were not able to compute the \(q\)-connectivity in feasible time due to large number of simplices.}
	\label{table:BlueBrain}
\end{table}

\textbf{C. elegans.} We set \(q=4\) for best comparison to the Blue Brain. The path lengths in the C. elegans network are drastically longer and show more variability. This implies in the first hand that the directed flag complex is relatively dense beyond dimension 4 to connect simplices along paths. As shown in Figure \ref{fig:simplex_counts}, the simplex counts in C. elegans drop less than two orders of magnitude after dimension 4, whereas in Blue Brain the drop is four orders of magnitude. Even though the number of simplices in C. elegans is very small compared to Blue Brain, the path structure still shows much longer connectivities. This might be related to C. elegans having a dense core in its network \cite{C_elegans} that Blue Brain is missing, or to the overrepresentation of reciprocal connections as we detail below.

The very low values of the path fractions for C. elegans indicate that the simplices along paths are rather strongly connected. To investigate this further we looked at the actual paths and observed that many simplices along the paths were spanned by the same vertices but ordered differently, corresponding to reciprocal edges in the digraph. Interestingly, \(q\)-analysis and the path fraction is thus able to pick out a known fact about the C. elegans network: there is an overrepresentation of reciprocally connected 2- and 3-cliques \cite{C_elegans}. The \(q\)-analysis shows further how the corresponding simplices connect to form paths and clusters of simplices. As suggested in \cite{C_elegans}, this strong connectivity may have a functional role in the C. elegans network.

Similarly to Blue Brain, we also computed simplicial paths with respect to all \((q,0,q+1)\)-connections for \(q \in \{1,\dots,6\}\). The path lengths and corresponding path fractions are shown in Table \ref{table:celegans}, which further shows the long path lengths in C. elegans, apart from the second to last simplicial dimension 6. 
\begin{table}[h!]
	\centering
	\begin{tabular}{|c|c|c|c|c|c|} 
		\hline
		(1,0,2) & (2,0,3) & (3,0,4) & (4,0,5) & (5,0,6) & (6,0,7) \\
		\hline
		19 & 22 & 29 & 26 & 24 & 2 \\ 
		\hline
		0.64 & 0.59 & 0.48 & 0.42 & 0.29 & 1 \\ 
		\hline
	\end{tabular}
	\caption{Simplicial path lengths and corresponding path fractions in the C. elegans network. The (6,0,7)-path is just a face inclusion of a 6-simplex into a 7-simplex and the path fraction is 1.}
	\label{table:celegans}
\end{table}

\textbf{Allen mouse and Drosophila.} For Allen mouse we set \(q=7\). The structure of the network is also characterised by quite long path lengths but less variability as compared to C. elegans. The same is true for the path fractions with more uniform distribution over the various connectivities. The actual network sizes of Allen mouse and Blue Brain are comparable, as well as the number of simplices used in the respective \(q\)-connectivity computations (recall that Blue Brain was restricted to Layers 1-4). Also in both cases the dimension of the flag complex is \(q+2\). It is then interesting that the Allen mouse, similarly to C. elegans, deviates so much from the more homogeneous structure of the Blue Brain.

For Drosophila we set \(q=8\). The path structure has very distinctive features with path lengths ranging from 3 to 44 and 53 at (8,0,5)- and (8,5,0)-connectivities, respectively. Similar variability is shown by the path fractions. 

The longest Drosophila paths at (8,0,5)- and (8,5,0)-connectivities both contain only 11 different vertices. The simplices along the paths all consist of those vertices but spanned by different reciprocal edges. The path analyses thus seem to converge on a structure resembling a completely connected directed 11-clique. Since we only used a subsample of the full Drosophila network this might indicate the existence of a high-dimensional (nearly) complete directed subgraph in the full network. Many of the paths with respect to other connectivities also consisted of the same 11 vertices spanning different simplices. The connectivities with path fraction 1 consist of vertices different from those 11.

\section*{Acknowledgements}
The author gratefully acknowledges Pedro Concei\c{c}\~{a}o, Joni Leino, and Barbara Mahler for providing feedback and corrections. Special thanks to Dejan Govc for a thorough reading and comments that polished many details in the exposition. Ran Levi provided support with many useful discussions. The author was funded by a collaboration agreement between the University of Aberdeen and EPFL, and by the KTH Royal Institute of Technology.
\bibliographystyle{plain}
\bibliography{bibliography}

\end{document}